\documentclass[reqno,12pt,draft%,dvips
]{amsart}
\usepackage{%showkeys,
amscd,amssymb,url}
\include{diagrams}

\allowdisplaybreaks[4]
\newarrow{TeXto}----{->}
\newcommand{\rto}{\rTeXto}

\newcommand{\rdto}{\rdTeXto}

\newcommand{\ldto}{\ldTeXto}

\renewcommand{\d}{\partial}
\newcommand{\g}{\mathfrak{g}}
\newcommand{\ep}{\varepsilon}

\newcommand{\fc}{\mathrm{fc}}
\newcommand{\idt}{\mathrm{id}}

\newcommand{\IL}{{\mathrm L}}
\newcommand{\cprime}{\/{\mathsurround=0pt$'$}}
\newcommand{\pinner}{\mathbin{\mathchoice
   {\hbox{\vrule width0.6em depth0pt height0.4pt
   \vrule width0.4pt depth0pt height0.8ex}}
   {\hbox{\vrule width0.6em depth0pt height0.4pt
   \vrule width0.4pt depth0pt height0.8ex}}
   {\hbox{\kern0.14em
   \vrule width0.48em depth0pt height0.4pt
   \vrule width0.4pt depth0pt height0.6ex\kern0.14em}}
   {\hbox{\kern0.1em
   \vrule width0.39em depth0pt height0.4pt
   \vrule width0.4pt depth0pt height0.5ex\kern0.1em}}}}
\newcommand{\inner}{\pinner\,}

\DeclareFontFamily{OML}{cyr}{}
\DeclareFontShape{OML}{cyr}{m}{n}{
   <5> <6> <7> <8> <9> gen * wncyr
   <10> <10.95> <12> <14.4> <17.28> <20.74> <24.88> wncyr10
  }{}
\DeclareSymbolFont{rusletters}{OML}{cyr}{m}{n}
\DeclareSymbolFontAlphabet{\rusmath}{rusletters}
\DeclareMathSymbol\re{\rusmath}{rusletters}{"03}

\DeclareMathOperator{\Lin}{Lin}

\DeclareMathOperator{\sym}{sym}

\DeclareMathOperator{\id}{id}
\DeclareMathOperator{\ad}{ad}

\newcommand{\bU}{\Bar{U}}

\newcommand{\BBR}{\mathbb{R}}

\newcommand{\CC}{\mathcal{C}}

\newcommand{\CE}{\mathcal{E}}
\newcommand{\CF}{\mathcal{F}}

\newcommand{\CU}{\mathcal{U}}
\newcommand{\CV}{\mathcal{V}}
\newcommand{\Dr}{\mathrm{D}}

\newcommand{\sbs}{\subset}
\newcommand{\wg}{\wedge}
\newcommand{\ot}{\otimes}

\newcommand{\la}{\lambda}
\newcommand{\Ga}{\Gamma}
\newcommand{\De}{\Delta}
\newcommand{\La}{\Lambda}

\newcommand{\Om}{\Omega}
\newcommand{\om}{\omega}
\newcommand{\si}{\sigma}
\newcommand{\vf}{\varphi}
\newcommand{\h}{-\hspace{0pt}}

\newcommand{\Ei}{\CE^\infty}
\newcommand{\Ci}{C^\infty}
\newcommand{\Ji}{J^\infty}

\def\ldb{{\rm [\![}}
\def\rdb{{\rm ]\!]}}
\newcommand{\fnij}[2]{\ldb{#1},{#2}\rdb}

\newcommand{\Gal}{\Ga_{\mathrm{loc}}}

\newtheorem{theorem}{Theorem}
\newtheorem{proposition}{Proposition}
\newtheorem{lemma}{Lemma}
\newtheorem{corollary}{Corollary}

\newtheorem*{theorem*}{Theorem}

\theoremstyle{definition}

\theoremstyle{remark}
\newtheorem{remark}{Remark}

\begin{document}
\title[On equations possessing flat representations]
{On symmetries and cohomological invariants\\
of equations possessing flat representations}

\author{S.~Igonin}
\address{S.~Igonin \\ Independent University of Moscow and
  University of Twente}
\curraddr{University of Twente \\ Department of Mathematics
\\ P.O.~Box 217 \\ 7500 AE Enschede \\  The Netherlands}
\email{igonin@mccme.ru}

\author{P.H.M.~Kersten}
\address{P.H.M.~Kersten \\ University of Twente \\ Department of Mathematics\\
P.O. Box 217 \\ 7500 AE Enschede \\ the Netherlands}
\email{kersten@math.utwente.nl}

\author{I.~Krasil{\cprime}shchik}
\address{I.~Krasil{\cprime}shchik \\ The Diffiety Institute, Moscow and
Independent University of Moscow\protect\newline
Correspondence to:
1st Tverskoy-Yamskoy per. 14, Apt. 45, 125047 Moscow, Russia}
\email{josephk@diffiety.ac.ru}
%\subjclass{37K35, 37K25, 58J10, 53C05}

\begin{abstract}
We study the equation $\CE_{\fc}$ of flat connections in a given fiber bundle
and discover a specific geometric structure on it, which we call a
\emph{flat representation}. We generalize this notion to arbitrary PDE and
prove that flat representations of an equation $\CE$ are in $1$-$1$
correspondence with morphisms $\vf\colon\CE\to\CE_{\fc}$, where $\CE$ and
$\CE_{\fc}$ are treated as submanifolds of infinite jet spaces. We show
that flat representations include several known types of zero\h curvature
formulations of PDEs. In particular, the Lax pairs of the self-dual
Yang--Mills  equations and their reductions are of this type. With each flat
representation $\vf$ we associate a complex $C_\vf$ of vector\h valued
differential forms such that $H^1(C_\vf)$ describes infinitesimal deformations
of the flat structure, which are responsible, in particular, for parameters
in B\"{a}cklund transformations. In addition, each higher infinitesimal
symmetry $S$ of $\CE$ defines a $1$\h cocycle $c_S$ of $C_\vf$. Symmetries with
exact $c_S$ form a subalgebra reflecting some geometric properties of $\CE$
and $\vf$. We show that the complex corresponding to $\CE_{\fc}$ itself is
$0$\h acyclic and $1$\h acyclic (independently of the bundle topology),
which means that higher symmetries of $\CE_{\fc}$ are exhausted by generalized
gauge ones, and compute the bracket on $0$\h cochains induced by commutation
of symmetries.
\end{abstract}
\maketitle
%%%%%%%%%%%%%%%%%%%%%%%%%%%%%%%%%%%%%%%%%%%%%%%%%

\tableofcontents

\section*{Introduction}
Nonlinear PDE in more than two independent variables with nontrivial invariance
properties is a sort of challenge for those who are interested in geometry of
differential equations and related topics. Evidently, chosen ``at random'', such
an equation will be of no interest, but \emph{natural} (whatever it means)
equations arising in physics and geometry may provide the researcher with the
desired results. Thinking of the second source (geometry), a first idea
is to consider equations describing geometrical structures (e.g., integrable
distributions, complex structures, etc.) and, by definition, possessing a rich
symmetry group. One of the simplest equations of this type is the equation
$\CE_{\fc}$ describing flat connections in a certain locally trivial fiber bundle.

Let $\pi\colon E\to M$, $\dim M=n$, $\dim E=n+m$, be a smooth locally trivial
fiber bundle over a smooth manifold $M$. A \emph{connection} $\nabla$ in the
bundle $\pi$ is a $\Ci(M)$\h linear correspondence that takes any vector field
$X$ on $M$ to a vector field $\nabla X$ on $E$ in such a way that
$\pi_*(\nabla X)_\theta=X_{\pi(\theta)}$ for any point $\theta\in E$.
A connection is called \emph{flat} if
\[
\nabla[X,Y]=[\nabla X,\nabla Y]
\]
for any two vector fields $X$ and $Y$ on $M$.

Let $\CU\sbs M$ be a local chart, such that $\pi$ becomes trivial over $\CU$,
with local coordinates $x_1,\dots,x_n$, $v^1,\dots,v^m$ in $\pi^{-1}(\CU)$,
$x_1,\dots,x_n$ being coordinates in $\CU$. Then $\nabla$ is determined by
$nm$ arbitrary functions $v_i^\alpha = v_i^\alpha(x_1,\dots,x_n,v^1,\dots,v^m)$
by the relations
\[
\nabla(\d x_i)=\d x_i+\sum_{\alpha=1}^mv_i^\alpha\d v^\alpha,\quad i=1,\dots,n,
\]
where and below $\d x_i=\d/\d x_i$, $\d v^\alpha=\d/\d v^\alpha$,
etc., while the condition of flatness is expressed in the form
\begin{multline*}
[\nabla(\d x_i),\nabla(\d x_j)]\equiv
\sum_{\alpha=1}^m\left(\frac{\d v_j^\alpha}{\d x_i}
+\sum_{\beta=1}^mv_i^\beta
\frac{\d v_j^\alpha}{\d v^\beta}\right)\d v^\alpha\\
-\sum_{\alpha=1}^m\left(\frac{\d v_i^\alpha}{\d x_j}
+\sum_{\beta=1}^mv_j^\beta
\frac{\d v_i^\alpha}{\d v^\beta}\right)\d v^\alpha=0,
\end{multline*}
$i,j=1,\dots,n$. Thus, flat connections are distinguished by the following system
of nonlinear equations
\begin{equation}\label{eq:fce}
\frac{\d v_j^\alpha}{\d x_i}+
\sum_{\beta=1}^mv_i^\beta\frac{\d v_j^\alpha}{\d v^\beta}
=\frac{\d v_i^\alpha}{\d x_j}+
\sum_{\beta=1}^mv_j^\beta\frac{\d v_i^\alpha}{\d v^\beta},
\end{equation}
$1\leq i<j\leq n$, $\alpha=1,\dots,m$.
This is a system of $mn(n-1)/2$ equations imposed on $mn$ unknown functions $v_i^\alpha$.

A strong motivation to study this equation is the following. If one assumes that
the coefficients $v_i^\alpha$ of a connection are expressed in terms of some
functions $f_k$ in a special way, then the flatness is equivalent to a system $\CE$
of differential equations on $f_k$ such that each solution to $\CE$ gives a flat
connection. It is well known that many  integrable nonlinear PDEs arise exactly
in this way. Treating the equations $\CE$ and $\CE_{\fc}$ with their differential
consequences as (infinite\h dimensional) submanifolds of the corresponding infinite
jet spaces, we obtain a smooth map
\begin{equation}\label{zcf}
  \vf\colon\CE\to\CE_{\fc},
\end{equation}
which preserves the Cartan distribution. Such a map generates on $\CE$ an additional
geometric structure which we call a \emph{flat representation}. In particular, the
identical map $\id\colon\CE_\fc\to\CE_\fc$ determines a canonical flat representation
for $\CE_\fc$ itself.

Thus $\CE_{\fc}$ is a sort of universal space where integrable equations live (or at least
a broad class of them). In particular, we show that coverings (generalized Lax pairs)
of equations in two independent variables as well as the Lax pairs of the self-dual
Yang--Mills equations and their reductions lead to flat representations. Therefore,
studying the geometry of $\CE_{\fc}$ may help to obtain new knowledge on integrable
systems, and in this paper we begin the study. Despite the fact that locally any flat
connection can be obtained from a trivial one by means of a local automorphism of the
bundle, there is no natural way to parameterize all flat connections, and this makes
the geometry of $\CE_{\fc}$ nontrivial.

We show that each flat representation~\eqref{zcf} leads to several complexes of
vector\h valued differential forms on $\CE$ constructed by means of the Nijenhuis
bracket. One subcomplex denoted by $C_\vf$ is of particular interest, since its
$1$\h cocycles are infinitesimal deformations of the flat representation, the exact
cocycles being trivial deformations. As usual, higher cohomology contains obstructions
to prolongation of infinitesimal deformations to formal deformations~\cite{gen-deform}.
In particular, zero\h curvature representations and B\"acklund transformations dependent
on a parameter provide an example of a deformed flat representation and determine,
therefore, some $1$\h cohomology classes of $C_\vf$.

In addition, each higher infinitesimal symmetry $S$ of $\CE$ defines a formal deformation
and, therefore, a $1$\h cocycle $c_S$ in $C_\vf$. Symmetries with exact cocycles $c_S$
form a Lie subalgebra with respect to the commutator. If morphism~\eqref{zcf} is an embedding
then it is precisely the subalgebra of those symmetries which are restrictions of some
symmetries of the equation of flat connections. If~\eqref{zcf} is coming from a
covering of $\CE$ then this subalgebra consists of symmetries which can be lifted to
the covering equation.

The complex $C_{\id}$ corresponding to the identical flat representation
$\id\colon\CE_{\fc}\to\CE_{\fc}$ itself is a new canonical complex associated with the bundle
$\pi\colon E\to M$. Symmetries of $\CE_{\fc}$ are in one-to-one correspondence with $1$\h cocycles
of this complex.

We prove that
\begin{equation}\label{h=0}
H^0(C_{\id})=H^1(C_{\id})=0
\end{equation}
independently of the topology of $M$ and $E$. This implies that the higher
symmetries of $\CE_{\fc}$ are in one-to-one correspondence with $0$\h cochains,
which are sections of some vector bundle over $\CE_{\fc}$ of rank $m$. That is,
symmetries are determined locally by $m$\h tuples of smooth functions on $\CE_{\fc}$.
This is similar to the case of an infinite jet space. Due to their coordinate expression,
it is natural to call them \emph{generalized gauge symmetries}. The commutator of
symmetries induces a bracket on $0$\h cochains, and we write down this bracket in
explicit terms.

Note that only for a few multidimensional equations the complete description of
symmetries is known. We mention the paper~\cite{anderson}, which proves that symmetries
of the vacuum Einstein equations in four spacetime dimensions are also freely determined
by $5$\h tuples of smooth functions on the equation.

The paper is organized as follows. In Section~\ref{sec:basics}, we review some facts
on the geometry of PDE and the Nijenhuis bracket. In Section~\ref{sec:eq-fc}, we
introduce the equation of flat connections and complexes of vector\h valued forms on it.
We prove~\eqref{h=0} and deduce the description of symmetries of $\CE_\fc$ in
Section~\ref{subs:result}. In Section~\ref{sec:flat-rep}, we introduce flat representations
in general and prove the above facts on their deformations and relations with symmetries.
The bracket on $0$\h cochains of $C_{\id}$ is computed in Subsection~\ref{subs:Lie-str}.
Finally, in Section~\ref{efr}, we show that coverings of PDE and the Lax pair of the
self-dual Yang--Mills equations lead to flat representations $\varphi$ and study the map
from symmetries to $1$\h cocycles of $C_\varphi$ for these examples.

\section{Basic facts}\label{sec:basics}
This is a brief review of the geometry of PDE. We refer to~\cite{KV-book,KK-book}
for more details.

\subsection{Notation agreements}\label{subs:notation}
Let $M$ be a smooth manifold and $\pi$ be a smooth fiber bundle over $M$. Throughout the
paper we use the following notation:
\begin{itemize}
\item
$\Dr(M)$ denotes the $\Ci(M)$\h module (and Lie algebra over $\BBR$) of vector fields
on $M$. A \emph{distribution} on $M$ is a $\Ci(M)$\h submodule of $\Dr(M)$.
\item
More generally, if $P$ is an arbitrary $\Ci(M)$\h module, then $\Dr(P)$ denotes the module
of $P$\h valued derivations:
\begin{multline*}
\Dr(P)=\{\,X\in\mathrm{Hom}_{\mathbb{R}}(\Ci(M),P)\mid \\
X(fg)=fX(g)+gX(f),\ f,g\in\Ci(M))\,\}.
\end{multline*}
\item
$\La^i(M)$ is the module of $i$\h differential forms on $M$ while
\[
\La^*(M)=\bigoplus_{i=1}^{\dim M}\La^i(M)
\]
is the corresponding algebra (with respect to the wedge product~$\wg$).
\item
$\Ga(\pi)$ (resp., $\Gal(\pi)$) is the set of all (resp., all local) sections of~$\pi$.
In particular, when $\pi$ is a vector bundle, $\Ga(\pi)$ is a $\Ci(M)$\h module.
\end{itemize}
\subsection{Jets and equations}\label{subs:jets}
Let $\pi\colon E\to M$ be a fiber bundle and $\theta\in E$, $\pi(\theta)=x\in M$. Consider
a local section $f\in\Gal(\pi)$ whose graph passes through the point $\theta$ and the class
$[f]_x^k$ of all local sections whose graphs are tangent to the graph of $f$ at $\theta$ with
order $\geq k$. The set
\[
J^k(\pi)=\{\,[f]_x^k\mid f\in\Gal(\pi),\ x\in M\,\}
\]
carries a natural structure of a smooth manifold and is called the \emph{manifold of $k$-jets}
of the bundle $\pi$. Moreover, the natural projections
\[
\pi_k\colon J^k(\pi)\to M,\ [f]_x^k\mapsto x,\quad
\pi_{k,k-1}\colon J^k(\pi)\to J^{k-1}(\pi),\ [f]_x^k\mapsto [f]_x^{k-1},
\]
are locally trivial bundles. If $\pi$ is a vector bundle, then $\pi_k$ is a vector bundle
as well while all $\pi_{k,k-1}$ carry a natural affine structure. So, we have an infinite
sequence of bundles
\[
\dots\to J^k(\pi)\xrightarrow{\pi_{k,k-1}}
J^{k-1}(\pi)\to\dots\to J^1(\pi)\xrightarrow{\pi_1}J^0(\pi)=E\xrightarrow{\pi}M
\]
and its inverse limit is called the \emph{manifold of infinite jets} and denoted by
$\Ji(\pi)$. In an obvious way, one can consider the fiber bundles
$\pi_\infty\colon\Ji(\pi)\to M$ and $\pi_{\infty,k}\colon\Ji(\pi)\to J^k(\pi)$. Using
the notation similar to the above one, we can state that $\Ji(\pi)$ consists of classes
of the form~$[f]_x^\infty$.

Let $f\in\Gal(\pi)$; then one can consider the sections $j_k(f)\in\Gal(\pi_k)$,
$k=0,1,\dots,\infty$, defined by $j_k(f)\colon x\mapsto [f]_x^k$ and called the $k$-jet
of $f$. Consider a point $\theta = [f]_x^\infty\in\Ji(\pi)$ and the tangent plane to the
graph of $j_\infty(f)$ at $\theta$. One can show that this plane does not depend on the
choice of $f$ and is determined by the point $\theta$ only. Denote this plane by~$\CC_\theta$.

Let $X\in\Dr(M)$. Then by the above construction there exists a unique vector field
$\CC X\in\Dr(\Ji(\pi))$ such that $\CC X_\theta\in\CC_\theta$ and $(\pi_\infty)_*\CC X_\theta
=X_x$ for any $\theta\in\Ji(\pi)$ and $x=\pi_\infty(\theta)$. Thus we obtain a connection
$\CC$ in the bundle $\pi_\infty$ called the \emph{Cartan connection}. The distribution
$\CC\Dr\sbs\Dr(\Ji(\pi))$ generated by the fields of the form $\CC X$, $X\in \Dr(M)$,
is called the \emph{Cartan distribution}.

\begin{proposition}\label{prop:Cdistr}
The Cartan distribution is integrable in the formal Frobenius sense\textup{,} i.e.\textup{,}
$[\CC\Dr,\CC\Dr]\sbs\CC\Dr$. The Cartan connection is flat\textup{,} i.e.\textup{,}
$[\CC X,\CC Y]=\CC[X,Y]$ for any $X,Y\in\Dr(M)$.
\end{proposition}

\begin{remark}\label{rem:Cdistr}
Consider a point $\theta_k\in J^k(\pi)$, $k>0$. Then it is completely determined by the
point $\theta_{k-1}=\pi_{k,k-1}(\theta_k)=[f]_{x}^{k-1}$, $f\in\Gal(\pi)$, and the tangent
plane to the graph of $j_{k-1}f$ at $\theta_{k-1}$. Consequently, sections of the bundle
$\pi_{k,k-1}\colon J^k(\pi)\to J^{k-1}(\pi)$ are identified with horizontal $n$\h dimensional
distributions on $J^{k-1}(\pi)$. In particular, sections of $\pi_{1,0}\colon J^1(\pi)\to E$
are just connections in the bundle~$\pi$.
\end{remark}

A $k$th order \emph{differential equation} imposed on sections of the bundle $\pi$ is a
submanifold $\CE^0\sbs J^k(\pi)$ of the manifold $J^k(\pi)$. Any equation may be represented
in the form
\[
\CE^0=\{\,\theta\in J^k(\pi)\mid (\De_\CE)_\theta=0,\ \De_\CE\in\Ga(\pi_k^*(\xi))\,\},
\]
where $\xi\colon E'\to M$ is a bundle and $\pi_k^*(\xi)$ is its pullback. The section
$\De_\CE=\De$ can be understood as a \emph{nonlinear differential operator} acting from
$\Ga(\pi)$ to $\Ga(\xi)$: $\De(f)=j_k(f)^*(\De_\CE)\in\Ga(\xi)$, $f\in\Ga(\pi)$. Of course,
neither $\xi$, nor $\De$ is unique.

The $l$th \emph{prolongation} of $\CE^0$ is the set
\begin{multline*}
\CE^l=\{\,[f]_x^{k+l}\mid\text{the graph of }j_k(f)
\text{ is tangent to }\CE^0\\
\text{ with order }\geq l\text{ at }[f]_x^k\in\CE\,\},
\end{multline*}
$l=0,1,\dots,\infty$. Restricting the maps $\pi_{k+l}$ and $\pi_{k+l,k+l-1}$ to $\CE^l$
and preserving the same notation for these restrictions, we obtain the sequence of maps
\[
\dots\to\CE^l\xrightarrow{\pi_{k+l,k+l-1}}
\CE^{l-1}\to\dots\to\CE^0\xrightarrow{\pi_{k,0}}E\xrightarrow{\pi}M.
\]
Without loss of generality, one can always assume $\pi_k$ to be surjective. Imposing natural
conditions of regularity, we shall also assume that all $\CE^l$ are smooth manifolds, while
the maps $\pi_{k+l,k+l-1}\colon \CE^l\to\CE^{l-1}$ are smooth fiber bundles (such equations
are called \emph{formally integrable}). Everywhere below the \emph{infinite prolongation}
$\Ei$ will be denoted by~$\CE$.

It can be shown that any vector field of the form $\CC X$, $X\in\Dr(M)$, is tangent to $\CE$.
This means that the Cartan connection restricts from the bundle $\Ji(\pi)\to M$ to $\CE\to M$
and determines on $\CE$ a horizontal $n$\h dimensional integrable distribution, the
\emph{Cartan distribution} on $\CE$. We denote by  $\CC\Dr_\CE\sbs\Dr(\CE)$ the submodule
generated by the vector fields lying in the Cartan distribution on $\CE$. Its maximal integral
manifolds are the infinite jets of solutions of $\CE$. A (\emph{higher}) \emph{symmetry} of
$\CE$ is a $\pi$\h vertical vector field on $\CE$ preserving the Cartan distribution, i.e.,
$S\in\Dr(\CE)$ is a symmetry if and only if
\begin{enumerate}
\item[(a)]
$S(f)=0$ for any function $f\in\Ci(M)\sbs\Ci(\CE)$;
\item[(b)]
$[S,Y]$ lies in $\CC\Dr_\CE$ whenever $Y\in\Dr(\CE)$ lies in $\CC\Dr_\CE$ (in fact, from (a)
it follows that $[S,\CC X]=0$ for any $X\in\Dr(M)$).
\end{enumerate}
The set $\sym\CE$ of all symmetries forms a Lie $\mathbb{R}$\h algebra with respect to the
commutation of vector fields. In the trivial case $\CE=\Ji(\pi)$ we denote this algebra by
$\sym\pi$. A complete description of symmetries is given by the following

\begin{theorem}\label{thm:sym-gen}
Let $\pi\colon E\to M$ be a locally trivial fiber bundle. Then there exists a one-to-one
correspondence between $\sym\pi$ and sections of the pullback $\pi_{\infty,0}^*(\pi^v)$\textup{,}
where $\pi^v\colon T^vE\to E$ is the vertical tangent bundle to $E$. If $\CE^0\sbs J^k(\pi)$
is an equation given by a differential operator $\De\colon\Ga(\pi)\to\Ga(\xi)$ and satisfying
the above formulated conditions\textup{,} then $\sym\CE$ is in one-to-one correspondence with
solutions of the equation
\begin{equation}\label{eq:sym-gen}
\ell_\CE(\vf)=0,
\end{equation}
where $\vf$ is a section of the pullback $\pi_{\infty,0}^*(\pi^v)$ restricted to $\CE$ while
$\ell_\CE$ is the linearization of $\De$ also restricted to $\CE$.
\end{theorem}

\begin{remark}\label{rem:gen-sect}
We call $\vf$ the \emph{generating section} of a symmetry. It may be considered as an element
of the module $\Ga(\pi^v)\otimes_{\Ci(E)}\Ci(\Ji(\pi))$. When $\pi$ is a vector bundle, this
module is isomorphic to $\Ga(\pi_\infty^*(\pi))=\Ga(\pi)\otimes_{\Ci(M)}\Ci(\Ji(\pi))$.
\end{remark}
\subsubsection*{Coordinates}
Let $\CU\sbs M$ be a local chart with coordinates $x_1,\dots,x_n$ and $\CU\times\mathbb{R}^m$ be
a trivialization of $\pi$ over $\CU$ with coordinates $v^1,\dots,v^m$ in $\mathbb{R}^m$. Then in
$\pi_\infty^{-1}(\CU)$ the so-called \emph{adapted coordinates} arise defined by
\begin{equation}\label{eq:adapted}
v_\sigma^\alpha([f]_x^l)=
\left.\frac{\d^l f^\alpha}{\d x_{i_1}\dots\d x_{i_l}}\right|_x,\quad
f\in\Gal(\pi),\ x\in\CU,\ \alpha=1,\dots,m,
\end{equation}
where $\sigma=i_1\dots i_l$ is a multi-index, $1\leq i_\alpha\leq n$. If $|\sigma|=l\leq k$, then
$v_\sigma^\alpha$ together with $x_1,\dots,x_n$ constitute a local coordinate system in
$\pi_k^{-1}(\CU)\sbs J^k(\pi)$.

In these coordinates, the Cartan connection is completely determined by its values on partial
derivatives and the corresponding vector fields are the so-called \emph{total derivatives}
\begin{equation}\label{eq:tot-der}
\CC(\d x_i)\equiv D_{x_i}=\d x_i + \sum_{\sigma,\alpha}v_{\sigma i}^\alpha\,\d v_\sigma^\alpha,
\end{equation}
$i=1,\dots,n$. The Cartan distribution is given by the system of \emph{Cartan forms}
\begin{equation}\label{eq:Cart-forms}
\om_\sigma^\alpha=dv_\sigma^\alpha-\sum_{i=1}^nv_{\sigma i}^\alpha\,dx_i,
\end{equation}
while the correspondence between $\Ga(\pi_\infty^*(\pi))$ and $\sym\pi$ (see
Theorem~\ref{thm:sym-gen} and Remark~\ref{rem:gen-sect}) is expressed by the formula
\begin{equation}\label{eq:evol-der}
\vf=(\vf^1,\dots,\vf^m)\mapsto\re_\vf=\sum_{\sigma,\alpha}D_\sigma(\vf^\alpha)\,
\d v_\sigma^\alpha,
\end{equation}
where $D_\sigma=D_{x_{i_1}}\circ\dots\circ D_{x_{i_k}}$ and $\re_\vf$ is called the
\emph{evolutionary vector field} with the \emph{generating section} $\vf$. Finally, if
$\De\colon\Ga(\pi)\to\Ga(\xi)$ is locally given by the relations
\[
u^s=F^s(x_1,\dots,x_n,\dots,v_\sigma^\alpha,\dots),\quad s=1,\dots,m',
\]
where $m'=\dim\xi$ and $u^s$ are coordinates in the fibers of the bundle $\xi$ over $\CU$, then
\begin{equation}\label{eq:lineariz}
\ell_\De=\left\Vert
\sum_\sigma\frac{\d F^s}{\d v_\sigma^\alpha}D_\sigma\right\Vert,\quad
\alpha=1,\dots,m,\ s=1,\dots,m',
\end{equation}
is its linearization. To get the expression for $\ell_\CE$ it suffices to
rewrite~\eqref{eq:lineariz} in \emph{internal coordinates} on $\CE$.

\subsection{Nijenhuis bracket}\label{subs:Nijenhuis}
Let $N$ be a smooth manifold and consider the $\mathbb{Z}$\h graded module $\Dr(\La^*(N))=
\oplus_{i=0}^{\dim N}\Dr(\La^i(N))$ over $\Ci(N)$ of $\La^*(N)$\h valued derivations. First
note that for any vector field $X\in\Dr(N)$ and a derivation $\Om\colon\Ci(N)\to\La^i(N)$
one can define a derivation $X\inner\Om\colon\Ci(N)\to\La^{i-1}(N)$ by
\[
(X\inner\Om)(f)=X\inner\Om(f),\quad f\in\Ci(N),
\]
(the inner product, or contraction), where $X\inner\Om(f)$ is the usual inner product of
a vector field and a differential form. Further, if $\om\in\La^j(N)$, then the inner product
$\Om\inner\om\in\La^{i+j-1}(N)$ is defined by induction on $i+j$:
\[
X\inner(\Om\inner\om)=(X\inner\Om)\inner\om-(-1)^i\Om\inner(X\inner\om)
\]
with an obvious base. Note that $\Om\inner$ is an operation of grading $i-1$.

Now, let us define the \emph{Lie derivative} $\IL_\Om\colon\La^j(N)\to\La^{j+i}(N)$ by
\[
\IL_\Om(\om)=d(\Om\inner\om)+(-1)^i\Om\inner(d\om).
\]
\begin{proposition}\label{prop:Nijen}
For any elements $\Om\in\Dr(\La^i(N))$ and $\Om'\in\Dr(\La^{i'}(N))$ there exists a
uniquely defined element $\fnij{\Om}{\Om'}\in\Dr(\La^{i+i'}(N))$ satisfying the identity
\[
[\IL_\Om,\IL_{\Om'}]=\IL_{\fnij{\Om}{\Om'}},
\]
where $[\IL_\Om,\IL_{\Om'}]$ denotes the graded commutator.
\end{proposition}
The element $\fnij{\Om}{\Om'}$ is called the \emph{Nijenhuis} (or \emph{Fr\"olicher--Nijenhuis})
\emph{bracket} of $\Om$ and $\Om'$. The module $\Dr(\La^*(N))$ is a graded Lie algebra with
respect to this bracket, i.e.,
\begin{gather}
\fnij{\Om}{\Om'}+(-1)^{ii'}\fnij{\Om'}{\Om}=0,\label{eq:skew}\\
\fnij{\Om}{\fnij{\Om'}{\Om''}}=\fnij{\fnij{\Om}{\Om'}}{\Om''}+
(-1)^{ii'}\fnij{\Om'}{\fnij{\Om}{\Om''}},\label{eq:Jacobi}
\end{gather}
where $i$, $i'$ and $i''$ are the degrees of $\Om$, $\Om'$ and $\Om''$, respectively. Note
that for ``decomposable" elements of $\Dr(\La^*(N))$, i.e., for elements of the form
$\Om=\om\otimes X$, $\om\in\La^*(N)$, $X\in\Dr(N)$, acting on $f\in\Ci(M)$ by $(\om\otimes X)(f)
=X(f)\om$, one has
\begin{multline}\label{eq:Nij-decomp}
\fnij{\om\otimes X}{\om'\otimes X'}=
\om\wg\om'\ot[X,X']+\om\wg X(\om')\ot X'-X'(\om)\wg\om'\ot X\\
+(-1)^i\,d\om\wg(X\inner\om')\ot X'+(-1)^i(X'\inner\om)\wg\,d\om'\ot X,
\end{multline}
where $X,X'\in\Dr(N)$, $\om\in\La^i(N)$, $\om'\in\La^{i'}(N)$.

Let $U\in\Dr(\La^1(N))$. Then the operator $\d_U=\fnij{U}{\cdot}\colon\Dr(\La^i(N))
\to\Dr(\La^{i+1}(N))$ is defined. If $\fnij{U}{U}=0$, then, due to properties~\eqref{eq:skew}
and~\eqref{eq:Jacobi}, the operator $\d_U$ satisfies the condition $\d_U\circ\d_U=0$
and thus we obtain the complex
\begin{equation}\label{eq:U-compl}
0\to\Dr(N)\xrightarrow{\d_U}\Dr(\La^1(N))\to\dots
\to\Dr(\La^i(N))\xrightarrow{\d_U}\Dr(\La^{i+1}(N))\to\dotsb
\end{equation}
with the corresponding cohomology groups $H_U^i(N)$.

\subsection{Coverings}\label{subs:coverings}
Let $\xi\colon N\to M$ and  $\xi'\colon N'\to M$ be two fiber bundles endowed with flat
connections $\nabla$ and $\nabla'$ respectively. Assume that $N'$ is fibered over $N$ by
some smooth surjection $\tau\colon N'\to N$ and $\xi'=\xi\circ\tau$. We say that the pair
$(N',\nabla')$ \emph{covers} the pair $(N',\nabla')$ (or $\tau$ is a \emph{covering}) if
for any $x\in M$, $X\in\Dr(M)$ and $\theta'\in N'$ one has $\tau_*(\nabla'X)_{\theta'}=
(\nabla X)_{\tau(\theta')}$. A covering $\tau$ is called \emph{linear} if
\begin{enumerate}
    \item $\tau\colon N'\to N$ is a vector bundle;
    \item vector fields of the form $\CC'X$, preserve the space $\Lin(\tau)\sbs\Ci(N')$, where
          $\Lin(\tau)$ consists of smooth functions on $N'$ linear along the fibers of $\tau$.
\end{enumerate}

If $N=\CE$, $\xi=\pi_\infty\colon\CE\to M$ and $\nabla=\CC$ is the Cartan connection on $\CE$,
we speak of a \emph{covering over $\CE$}. In local coordinates, any covering structure over
$\CE$ is determined by a set of $\tau$\h vertical vector fields $X_1,\dots,X_n$ on $N'$ such
that the conditions
\[
[D_{x_i}+X_i,D_{x_j}+X_j]=0
\]
hold for all $1\le i<j\le n$, where $D_{x_\alpha}$ are the total derivatives given
by~\eqref{eq:tot-der}. In this case, if $\tau$ is finite\h dimensional then the manifold $N'$ is
locally also isomorphic to an infinitely prolonged differential equation called the
\emph{covering equation} with the Cartan connection $\nabla'$ and the total derivatives
$D_{x_i}+X_i$.

A detailed discussion of the theory of coverings over nonlinear PDE can be found
in~\cite{KV-book, KV-trends}.

\subsection{Horizontal cohomologies}\label{subs:hor-coh}
The notion of the Cartan connection can be generalized. Namely, consider an infinitely prolonged
equation $\pi_\infty\colon\CE\to M$ and two vector bundles $\xi$, $\xi'$ over $M$. Let
$\De\colon\Ga(\xi)\to\Ga(\xi')$ be a linear differential operator. Consider the pullbacks
$\pi_\infty^*(\xi)$, $\pi_\infty^*(\xi')$ and denote the corresponding modules of sections by
$\Ga(\pi,\xi)$, $\Ga(\pi,\xi')$, respectively. Then there exists a unique linear differential
operator $\CC\De\colon\Ga(\pi,\xi)\to\Ga(\pi,\xi')$ satisfying the local condition
\[
j_\infty(f)^*\CC\De(\vf)=\De(j_\infty(f)^*(\vf))
\]
for any $f\in\Gal(\pi)$ such that $\theta=[f]_x^\infty\in\CE$ and any $\vf\in\Gal(\pi,\xi)$ in
the neighborhood of~$\theta$.

In particular, if $\xi=\xi'\colon\bigwedge^iT^*M\to M$ is the $i$th external power of the
cotangent bundle, elements of the module $\Ga(\pi,\xi)\sbs\La^i(\CE)$ are called \emph{horizontal
forms} while the module itself is denoted by $\La_h^i(\CE)$. Taking the de~Rham differential
$d\colon\La^i(M)\to\La^{i+1}(M)$ and introducing the notation $\CC d=d_h$, we obtain the complex
\[
0\to\Ci(\CE)\xrightarrow{d_h}\La_h^1(\CE)\to\dots\to\La_h^i(\CE)\xrightarrow{d_h}\La_h^{i+1}(\CE)
\to\dots\to\La_h^n(\CE),
\]
the \emph{horizontal de~Rham complex} of $\CE$. Its cohomology at the $i$th term is denoted by
$H_h^i(\CE)$ and called the \emph{horizontal cohomology}. In the adapted coordinates we have
\begin{equation}\label{d_h}
d_h(fdx_{i_1}\wg\dots\wg\,dx_{i_k})=\sum_{i=1}^nD_{x_i}(f)\,dx_i\wg\,dx_{i_1}\wg\dots\wg\, dx_{i_k},
\quad f\in C^\infty(\CE).
\end{equation}

\begin{proposition}\label{prop:hor-triv}
If $\CE=\Ji(\pi)$\textup{,} $\pi\colon E\to M$\textup{,} then $H_h^i(\CE)$ for $i<n$ is isomorphic
to the $i$th de Rham cohomology of the manifold $E$.
\end{proposition}
\begin{remark}\label{rem:hor-coeff}
Let $\tau$ be a linear covering over $\CE$ and $P=\Ga(\tau)$. Then one can extend the horizontal
differential to $d_h\colon\La_h^i(\CE)\ot P\to\La_h^{i+1}(\CE)\ot P$ and thus obtain horizontal
cohomologies with coefficients in~$P$.
\end{remark}
For more details on horizontal cohomology see~\cite{KrasVer}.

\section{The equation of flat connections}\label{sec:eq-fc}

\subsection{More on connections}\label{subs:more-conn}
Let $\pi\colon E\to M$ be a fiber bundle. Let us consider the $\Ci(E)$\h module $\Dr(M,E)$ of
$\Ci(E)$\h valued derivations $\Ci(M)\to\Ci(E)$. Due to the embedding $\pi^*\colon\Ci(M)\to\Ci(E)$,
for any vector field $X\in\Dr(E)$ one can consider its restriction $X_M\in\Dr(M,E)$ to $M$. A
\emph{connection} in the bundle $\pi$ is a $\Ci(M)$\h linear map $\nabla\colon\Dr(M)\to\Dr(E)$
such that $\nabla(X)_M=X$ for any $X\in\Dr(M)$.

\begin{remark}\label{def:conn}
If $M$ is finite\h dimensional, this definition becomes equivalent to the following one: a
connection in a bundle $\pi$ is a $\Ci(E)$\h linear map $\nabla\colon\Dr(M,E)\to\Dr(E)$ such
that $\nabla(X)_M=X$ for any $X\in\Dr(M,E)$. In particular, we use this fact to define the
connection form below.
\end{remark}
Consider the element $\bU_\nabla\in\Dr(\La^1(E))$ defined by
\[
X\inner\bU_\nabla(f)=(\nabla X_M)(f),
\]
where $X\in\Dr(E)$, $f\in\Ci(E)$. The element $U_\nabla\in\Dr(\La^1(E))$ defined by
\[
X\inner U_\nabla(f)=X(f)-(\nabla X_M)(f)
\]
is called the \emph{connection form} of $\nabla$. Obviously, $\bU_\nabla+U_\nabla=d$, where
$d\colon\Ci(E)\to\La^1(E)$ is the de~Rham differential.

Recall that a connection $\nabla$ is called \emph{flat}, if $[\nabla X,\nabla Y]=\nabla[X,Y]$
for any $X,\,Y\in\Ci(M)$.
\begin{proposition}\label{prop:curvature}
For any $\Om\in\Dr(\La^*(E))$ one has $\fnij{U_\nabla}{\Om}=\fnij{\bU_\nabla}{\Om}$\textup{,}
and $\nabla$ is flat if and only if $\fnij{\bU_\nabla}{\bU_\nabla}=\fnij{U_\nabla}{U_\nabla}=0$.
\end{proposition}

{}From Proposition~\ref{prop:curvature} and Subsection~\ref{subs:Nijenhuis} it follows that to
any flat connection $\nabla$ we can associate the complex
\begin{equation}\label{eq:nabla-compl}
0\to\Dr(E)\xrightarrow{\d_\nabla}\Dr(\La^1(E))\to\dots
\to\Dr(\La^i(E))\xrightarrow{\d_\nabla}\Dr(\La^{i+1}(E))\to\dotsb,
\end{equation}
where
$\d_\nabla=\fnij{\bU_\nabla}{\cdot}=-\fnij{U_\nabla}{\cdot}$, cf.~\eqref{eq:U-compl}, and the
corresponding cohomology groups $H_\nabla^i(E)$.

\begin{remark}\label{rem:split}
In any fiber bundle $\pi\colon E\to M$ one can consider \emph{$\pi$\h vertical} vector fields
defined by
\[
\Dr^v(E)=\{\,X\in\Dr(E)\mid X_M=0\,\}.
\]
Dually, \emph{horizontal} $q$-forms can be introduced as
\[
\La_h^q(E)=\{\,\om\in\La^q(E)\mid X_1\inner\dots\inner X_q\inner\om=0,\
X_1,\dots,X_q\in\Dr^v(E)\,\}.
\]
Given a flat connection $\nabla$, one can split the module $\Dr(E)$ into the direct sum
$\Dr(E)=\Dr^v(E)\oplus\nabla\Dr(E)$, where $\nabla\Dr(E)\sbs\Dr(E)$ is the submodule generated
by the vector fields $\nabla X$, $X\in\Dr(M)$. This splitting induces the dual one:
$\La^1(E)=\La_h^1(E)\oplus\La_\nabla^1(E)$, where $\La_\nabla^1(E)$ is the annihilator of
$\nabla\Dr(E)$. Finally we obtain the splitting $\La^i(E)=\oplus_{i=p+q}\La^{p,q}(E)$, where
$\La^{p,q}(E)=\La_\nabla^p(E)\otimes\La_h^q(E)$ and $\La_\nabla^p(E)$ is defined in an obvious
way.
\begin{proposition}\label{prop:split}
Let $\d_\nabla$ be the differential from~\eqref{eq:nabla-compl}. Then
\[
\d_\nabla(\Dr(\La^{p,q}(E))\sbs\Dr(\La^{p,q+1}(E))
\]
and  thus the complexes
\begin{multline}\label{eq:nabla-compl-p}
0\to\Dr(\La^p(E))\xrightarrow{\d_\nabla}\Dr(\La^{p,1}(E))\to\dotsb\\
\dots\to\Dr(\La^{p,q}(E))
\xrightarrow{\d_\nabla}\Dr(\La^{p,q+1}(E))\to\dotsb
\end{multline}
are defined.
\end{proposition}
The corresponding cohomology is denoted by $H_\nabla^{p,q}(E)$.
\end{remark}

\begin{remark}\label{rem:subcompl}
In the sequel, we shall also deal with subcomplexes of complexes~\eqref{eq:nabla-compl}
and~\eqref{eq:nabla-compl-p}. As an example, let us indicate the complex
\begin{multline}\label{eq:nabla-compl-v}
0\to\Dr^v(E)\xrightarrow{\d_\nabla}\Dr^v(\La^1(E))\to\dotsb\\
\dots\to\Dr^v(\La^i(E))\xrightarrow{\d_\nabla}\Dr^v(\La^{i+1}(E))\to\dotsb
\end{multline}
of vertical $\La^*(E)$\h valued derivations.
\end{remark}

\subsubsection*{Local coordinates}
Let $\CU\sbs M$ be a local chart with the coordinates $x_1,\dots,x_n$, $v^1,\dots,v^m$ in
$\pi^{-1}(\CU)\sbs E$ (the case $m=\infty$ is included). Then, as it was noted above,
a connection $\nabla$ in $\pi$ is given by the correspondence
\[
\d x_i\mapsto\nabla_{x_i}=\nabla(\d x_i)=
\d x_i+\sum_{\alpha=1}^mv_i^\alpha\,\d v^\alpha,\quad
i=1,\dots,n,\ v_i^\alpha\in\Ci(E).
\]
By definition, one has
\begin{equation}\label{eq:coord-conn-form}
\bU_\nabla=\sum_{i=1}^ndx_i\otimes\nabla_{x_i},\quad
U_\nabla=
\sum_{\alpha=1}^m\Big(dv^\alpha-\sum_{i=1}^nv_i^\alpha\,dx_i\Big)\otimes\d v^\alpha.
\end{equation}
Restricting ourselves to vertical complex~\eqref{eq:nabla-compl-v}, from~\eqref{eq:Nij-decomp}
and~\eqref{eq:coord-conn-form} we get
\begin{equation}\label{eq:coord-differ}
\d_\nabla\Big(\sum_\alpha\theta^\alpha\otimes\d v^a\Big)
=\sum_{i=1}^ndx_i\wg\sum_{\alpha=1}^m\Big(\nabla_{x_i}(\theta^\alpha)
-\sum_{\beta=1}^m\frac{\d v_i^\alpha}
{\d v^\beta}\theta^\beta\Big)\otimes\d v^\alpha,
\end{equation}
where $\theta^\alpha\in\La^*(E)$.

\subsection{A coordinate-free formulation}\label{subs:coord-free}
Fix a finite\h dimensional bundle $\pi\colon E\to M$ and recall that connections (see
Remark~\ref{rem:Cdistr}) in this bundle are identified with the sections of the bundle
$\pi_{1,0}\colon J^1(\pi)\to E$. Consider a point $\theta_1\in J^1(\pi_{1,0})$. Then it can be
represented as a class $[\nabla]_\theta^1$, $\theta=(\pi_{1,0})(\theta_1)$, of local sections of
$\pi_{1,0}$, or local connections in the bundle $\pi$ over a neighborhood of the point
$\pi_1(\theta)\in M$.

Clearly, the value of $\fnij{\bU_\nabla}{\bU_\nabla}$ at $\theta$ is independent of the choice of
a representative $\nabla$ in the class $\theta_1=[\nabla]_\theta^1$ and is determined by the point
$\theta_1\in  J^1(\pi_{1,0})$ only. The \emph{equation of flat connections} is defined to be the
submanifold
\[
\CE=\{\theta_1\in J^1(\pi_{1,0})\mid
\fnij{\bU_\nabla}{\bU_\nabla}_{\pi_{1,0}(\theta_1)}=0,\ \theta_1=[\nabla]_\theta^1\}.
\]
Indeed, by Proposition~\ref{prop:curvature}, a section $s$ of $\pi_{1,0}$ determines a flat
connection if and only if $s(E)\subset\CE$.

In adapted coordinates this equation can be rewritten as
\begin{equation}\label{eq:fce-brack3}
\fnij{\sum_i dx_i\otimes\Big(D_{x_i}+\sum_\alpha v_i^\alpha D_{v^\alpha}\Big)}
{\sum_i dx_i\otimes\Big(D_{x_i}+\sum_\alpha v_i^\alpha D_{v^\alpha}\Big)}=0,
\end{equation}
or, equivalently,
\[
[D_{x_i}+\sum_\alpha v_i^\alpha D_{v^\alpha},\,
D_{x_j}+\sum_\alpha v_j^\alpha D_{v^\alpha}]=0,\quad 1\le i<j\le n,
\]
where $D_{x_i}$, $i=1,\dots,n$, $D_{v^\alpha}$, $\alpha=1,\dots,m$, are the total derivatives
corresponding to $\d x_i$ and $\d v^\alpha$, respectively, by the Cartan connection
$\CC$ in the fiber bundle $(\pi_{1,0})_\infty\colon\Ji(\pi_{1,0})\to E$.

We shall now describe the element $\sum_i\big(D_{x_i}+\sum_\alpha v_i^\alpha D_{v^\alpha}\big)
\otimes dx_i$ appearing in equation~\eqref{eq:fce-brack3} in a new, invariant way. To this end,
consider the Cartan connection $\CC^\pi$ in the bundle $\pi_\infty\colon\Ji(\pi)\to M$. This
connection takes any vector field $X\in M$ to the vector field $\CC^\pi X$ on $\Ji(\pi)$. Due
to the projection $\pi_{\infty,0}\colon\Ji(\pi)\to E$, we can restrict $\CC^\pi X$ to $E$
(more precisely, we restrict the map $\CC^\pi X\colon\Ci(\Ji(\pi))\to\Ci(\Ji(\pi))$ to the
subalgebra $\Ci(E)\sbs\Ci(\Ji(\pi))$). This restriction, $(\CC^\pi X)_0$, is an element of the
module $\Dr(E,J^1(\pi))$, and, by Remark~\ref{def:conn}, since $\dim E<\infty$, we can apply the
Cartan connection $\CC$ thus obtaining a vector field on $\Ji(\pi_{1,0})$. The sequence of maps
\[
\Dr(M)\xrightarrow{\CC^\pi}\Dr(\Ji(\pi))\xrightarrow{\text{restriction to } E}
\Dr(E,J^1(\pi))\xrightarrow{\CC}\Dr(\Ji(\pi_{1,0}))
\]
results in a new connection $\CC_{\fc}\colon\Dr(M)\to\Dr(\Ji(\pi_{1,0}))$
in the bundle
\[
\Tilde{\pi}=\pi\circ(\pi_{1,0})_\infty\colon\Ji(\pi_{1,0})\to M.
\]
When dealing with this connection, we shall use the notation $\bar{U}_{\fc}=
\bar{U}_{\CC_{\fc}}$ and $\d_{\fc}=\d_{\CC_{\fc}}$.

It is easily checked that in local coordinates one has
\[
\bU_{\fc}=\sum_i dx_i\otimes\big(D_{x_i}+\sum_\alpha v_i^\alpha D_{v^\alpha}\big)
\]
and we have the following
\begin{proposition}\label{prop:fce-new}
The equation of flat connections $\CE\sbs J^1(\pi_{1,0})$ is distinguished by the condition
$\fnij{\bU_{\fc}}{\bU_{\fc}}_{\theta_1}=0$\textup{,} $\theta_1\in J^1(\pi_{1,0})$.
All vector fields of the form $\CC_{\fc}X$\textup{,} $X\in\Dr(M)$\textup{,} are tangent to
its infinite prolongation $\CE_\fc\sbs\Ji(\pi_{1,0})$
and thus the restriction of $\CC_{\fc}$ to the subbundle $\tilde\pi\colon\CE_\fc\to M$ of
the bundle $\pi\circ(\pi_{1,0})_\infty\colon\Ji(\pi_{1,0})\to M$ is flat.
\end{proposition}

\section{Symmetries of the equation of flat connections}\label{subs:result}
{}From the last proposition it follows that the connection $\CC_{\fc}$ determines on $\CE_\fc$
the complexes similar to~\eqref{eq:nabla-compl-p}. We are especially interested in the case
$p=0$ and thus obtain
\begin{multline}\label{eq:CC-compl}
0\to\Dr(\CE_\fc)\xrightarrow{\d_{\fc}}\Dr(\La_h^1(\CE_\fc))\xrightarrow{\d_{\fc}}\dotsb\\
\dots\to\Dr(\La_h^q(\CE_\fc))\xrightarrow{\d_{\fc}}\Dr(\La_h^{q+1}(E))\to\dotsb,
\end{multline}
where $\d_\fc=\d_{\CC_\fc}$.

Consider now in $\CC\Dr(\CE_\fc)$ the submodule of $\tilde\pi$\h vertical vector fields.
In coordinates, it is spanned by the vector fields $D_{v^\alpha}$,
$\alpha=1,\dots,m$.
Set also $\CV^q=\tilde\pi^*(\La^q(M))\otimes\CV\sbs\Dr(\La_h^q(\CE_\fc))$.
\begin{lemma}\label{lem:subcompl}
The modules $\CV^q$ are preserved by the differential $\d_{\fc}$\textup{,} i.e.\textup{,}
$\d_{\fc}(\CV^q)\sbs\CV^{q+1}$.
\end{lemma}
\begin{proof}
This is easily seen from the definition of $\CV^q$. By~\eqref{eq:coord-differ}, in coordinates
$\d_{\fc}$ acts as follows
\begin{multline}\label{eq:CC-differ}
\d_{\fc}(f\,dx_{i_1}\wg\dots\wg\,dx_{i_q}\otimes D_{v^\alpha})=
\sum_{i=1}^n D_i(f)\,dx_i\wg\,dx_{i_1}\,\wg\dots\wg\,dx_{i_q}\otimes D_{v^\alpha}\\
-(-1)^q\sum_{i=1}^n
\sum_{\beta=1}^m D_{v^\alpha}
(v_i^\beta)f\,dx_i\,\wg dx_{i_1}\wg\dots\wg\,dx_{i_q}\otimes D_{v^\beta},
\end{multline}
where $f\in\Ci(\CE_\fc)$ and, as everywhere below, $D_i=D_{x_i}+\sum_\beta v_i^\beta D_{v^\beta}
=\CC_{\fc}\d x_i$ is the total derivative with respect to the connection $\CC_{\fc}$.
\end{proof}
Consequently, we obtain the complex
\begin{equation}\label{eq:CC-Vcompl}
0\to\CV^0=\CV\xrightarrow{\d_{\fc}}\CV^1\xrightarrow{\d_{\fc}}\dotsb\\
\dots\to\CV^q\xrightarrow{\d_{\fc}}\CV^{q+1}\to\dotsb
\end{equation}

Let us now return back to the equation of flat connections written in the
form~\eqref{eq:fce-brack3}. Since this equation is bilinear, its linearization in local
coordinates is of the form
\begin{multline*}
\fnij{\bU_{\fc}}{\sum_{i,\alpha}\vf_i^\alpha dx_i\otimes D_{v^\alpha}}
+\fnij{\sum_{i,\alpha}\vf_i^\alpha dx_i\otimes D_{v^\alpha}}{\bU_{\fc}}=\\
=2\fnij{\bU_{\fc}}{\sum_{i,\alpha}\vf_i^\alpha dx_i\otimes D_{v^\alpha}}=0,
\end{multline*}
where $\vf=\sum_{i,\alpha}\vf_i^\alpha dx_i\otimes D_{v^\alpha}$ is the generating section of
a symmetry, $\vf_i^\alpha\in\Ci(\CE_\fc)$, or
\begin{equation}\label{eq:fce-sym}
\ell_{\CE_\fc}\vf\equiv 2\d_{\fc}\vf=0,\quad\vf\in\CV^1.
\end{equation}

\begin{proposition}\label{prop:1-coc}
The map
\begin{equation}\label{sym-co}
        \sym\CE_\fc\to\CV^1,\quad S\mapsto\fnij{\bU_{\fc}}{S}
\end{equation}
is an isomorphism of the space of symmetries onto the space of $1$\h cocycles of
complex~\eqref{eq:CC-Vcompl}.
\end{proposition}
\begin{proof}
Let $S\in\sym \CE_\fc$. Since $S$ commutes with all total derivatives, we have
$\fnij{\bU_{\fc}}{S}\in\CV^1$. Moreover, $\fnij{\bU_{\fc}}{S}$ is a cocycle, since
$\fnij{\bU_{\fc}}{\fnij{\bU_{\fc}}{\cdot}}=0$ (since $S\notin\CV$, this cocycle may be not exact).
The fact that~\eqref{sym-co} is an isomorphism onto the kernel of $\d_{\fc}$ follows from the
above coordinate consideration, since $\fnij{\bU_{\fc}}{S}=-\sum_{i,\alpha}S(v_i^\alpha)
dx_i\otimes D_{v^\alpha}$ and $\vf_i^\alpha=S(v_i^\alpha)$ is nothing but the generating
section of the symmetry~$S$.
\end{proof}

\begin{remark}\label{rem:fce-genfunc}
It is not surprising that in~\eqref{eq:fce-sym} the symmetry generating section is an element
of the module $\CV^1$. Indeed, in our case the bundle $\xi=\pi_{1,0}\colon J^1(\pi)\to E$ is
affine and $\pi_{0,1}^v$ (see Remark~\ref{rem:gen-sect}) is isomorphic to the pullback
$\pi_1^*(\pi\otimes\tau^*)$, where $\tau^*\colon T^*M\to M$ is the cotangent bundle. Consequently,
the module of generating sections is isomorphic to $\CV^1$.
\end{remark}
\begin{theorem}\label{thm:main}
Complex~\eqref{eq:CC-Vcompl} is $0$\h acyclic\textup{,} i.e.\textup{,}
\begin{equation}
\label{0ac}
\ker\big(\d_{\fc}\colon\CV\to\CV^1\big)=0,
\end{equation}
and $1$\h acyclic\textup{,} i.e.\textup{,}
\begin{equation}
\label{1ac}
\d_{\fc}(\CV)=\ker\big(\d_{\fc}\colon\CV^1\to\CV^2\big).
\end{equation}
\end{theorem}
\begin{corollary}[description of symmetries]\label{cor:sym}
Any symmetry $(\vf_i^\alpha)$\textup{,} $i=1,\dots,n$\textup{,} $\alpha=1,\dots,m$\textup{,}
of equation~\eqref{eq:fce} is of the form
\begin{equation}\label{eq:spec-coord}
\vf_i^\alpha=D_i(f^\alpha)-\sum_{\beta=1}^m v_i^{\alpha,\beta}f^\beta
\end{equation}
for arbitrary functions $f^\alpha\in\Ci(\CE_\fc)$ uniquely determined by $(\vf_i^\alpha)$, where
$v_i^{\alpha,\beta}=D_{v^\beta}(v_i^{\alpha})$.
\end{corollary}
In Subsection~\ref{subs:Lie-str}, we shall describe a Lie algebra structure of symmetries in
terms of functions $f^\alpha$.
\begin{proof}[Proof of Corollary~\textup{\ref{cor:sym}}]
The result immediately follows from Theorem~\ref{thm:main} in combination
with~\eqref{eq:CC-differ} and~\eqref{eq:fce-sym}.
\end{proof}

To prove Theorem~\ref{thm:main}, we shall need some auxiliary constructions.
Consider the following smooth functions on $\CE_\fc$
\begin{equation}\label{eq:v-coord}
v_{I}^{\alpha,A}=D_{v^A}D_I(v^\alpha),\quad |I|>0,
\end{equation}
where $D_{v^A}=D_{v^{\alpha_1}}\circ\dots\circ D_{v^{\alpha_l}}$,
$D_I=D_{i_1}\circ\dots\circ D_{i_k}$ for the
multi\h indices
$A=\alpha_1\dots\alpha_l$
and $I=i_1\dots i_k$, respectively. In
particular, $v^{\alpha,\varnothing}_\varnothing=v^\alpha$ and $v_i^{\alpha,\varnothing}=v_i^\alpha$.
It is easily seen that the functions $x_i,\,v_I^{\alpha,A}$ form a system of coordinates on
$\CE_\fc$. Let $\CF_r\subset\Ci(\CE_\fc)$ be the subalgebra of functions dependent on the
variables
\[
x_i,\ v_{I}^{\alpha,A},\quad |I|\le r,
\]
only. Then we have
\[
 \CF_r\subset\CF_{r+1},\quad  \bigcup_{r\ge 0}\CF_r=\Ci(\CE_\fc).
\]
Evidently, for $f\in\CF_r$ one obtains
\begin{equation} \label{eq:diff-mod}
D_i(f)\equiv\sum_{\alpha,A,I}v_{Ii}^{\alpha,A}\frac{\d f}{\d v_{I}^{\alpha,A}}
\mod \CF_r.
\end{equation}
We have also the filtration $\CV^q_r\subset \CV^q_{r+1}$ in the modules $\CV^q$, where $\CV^q_r
\subset\CV^q$ consists of forms with coefficients in $\CF_r$. {}From~\eqref{eq:CC-differ}
and~\eqref{eq:diff-mod} for $r>0$ we get
\begin{multline} \label{eq:diff-mod1}
\d_{\fc}(f\,dx_{l_1}\wg\dots\wg\,dx_{l_q}\otimes D_{v^\alpha}) \\
\equiv  \sum_{i=1}^n\sum_{I,\beta,A} v_{Ii}^{\beta,A}
\frac{\d f}{\d v_{I}^{\beta,A}}\,dx_i\wg\,dx_{l_1}\wg\dots\wg\,dx_{l_q}
\otimes D_{v^\alpha}\mod \CV^{q+1}_r.
\end{multline}

On the other hand, consider the infinite\h dimensional vector bundle
\[
\pi'\colon\BBR^\infty\times\BBR^n\to\BBR^n
\]
with the coordinates $u^{\alpha,A}$,
$\alpha=1,\dots,m$, $A=\alpha_1\dots\alpha_l$ where $\alpha_j\ge 0$, along fibers and coordinates
$y_i$, $i=1,\dots,n$, in the base. The infinite jet bundle $\Ji(\pi')\to\BBR^n$ of $\pi'$ has the
adapted coordinates
\[
u^{\alpha,A}_{I}=D_{y_I}(u^{\alpha,A})
\]
along the fibers, where, as above,
$D_{y_I}$ denotes the composition $D_{y_{i_1}}\circ\dots\circ D_{y_{i_k}}$
of total derivatives.

To any element
\[
s=\sum_{i,\alpha}f^\alpha_i\,dx_i\otimes D_{v^\alpha}\in \CV^1
\]
and an integer $a\ge 0$ we associate an $m$-tuple $\om(s,a)^1,\dots,\om(s,a)^m$ of
parameter\h dependent horizontal $1$-forms on $\Ji(\pi')$ as follows.

Let $r(s)$ be the minimal integer such that all the coefficients $f^\alpha_i$ of $s$
belong to $\CF_{r(s)}$. Replace the variables $v_{I}^{\beta,A}$ with $|I|\ge r(s)-a$ in
$f^\alpha_i$ by $u_I^{\alpha,A}$ and denote thus obtained functions by $\tilde f^\alpha_i$.
We treat $\tilde f^\alpha_{i}$ as functions on $\Ji(\pi')$ dependent on the parameters
$x_1,\dots,x_n$ and $v_{I}^{\beta,A}$ with $|I|<r(s)-a$. Set
\[
 \om(s,a)^\alpha=\sum_i\tilde f^\alpha_{i}\,dy_i.
\]
We say that an element $s\in\CV^1$ is \emph{semilinear} if the coefficients of $s$ are linear
with respect to the highest order coordinates $v_{I}^{\alpha,A}$, $|I|=r(s)$.

\begin{lemma}\label{lem:d(w)}
The following facts are valid\textup{:}
\begin{enumerate}
\item\label{lem:d(w)-1}
If $s\in\CV^1$ and $\d_{\fc}(s)=0$ then $d_h(\om(s,0)^\alpha)=0$ for each $\alpha$ and for
all values of the parameters\textup{,} where $d_h$ is the horizontal de~Rham differential
on $J^\infty(\pi')$.
\item\label{lem:d(w)-2}
If $s\in\CV^1$ is semilinear\textup{,} $r(s)\ge 2$\textup{,}
and $\d_{\fc}(s)=0$\textup{,} then
\[
d_h(\om(s,1)^\alpha)=0.
\]
\end{enumerate}
\end{lemma}

\begin{proof}
(1)
By~\eqref{eq:diff-mod1} and~\eqref{d_h}, the equations $d_h(\om(s,0)^\alpha)=0$ for
all $\alpha=1,\dots,m$ are equivalent to
\[
\d_{\fc}(s)=0\mod \CV^2_{r(s)}.
\]

(2)
Let $\tilde\CF_r\subset\CF_r$ be the $\CF_{r-1}$\h submodule generated by
$v_{I}^{\alpha,A}$, $|I|=r$, and, respectively, let $\tilde\CV^1_r\subset\CV^1_r$
be the subset of forms with coefficients in $\tilde\CF_r$. Then from the identity
\begin{equation}\label{eq:commut}
[D_i,D_{v^\alpha}]=-\sum_\beta v_i^{\beta,\alpha}\,D_{v^\beta}
\end{equation}
we obtain
\begin{equation} \label{eq:D(C)}
D_i(v_{I}^{\alpha,A})\equiv v_{Ii}^{\alpha,A} \mod \tilde\CF_{|I|}.
\end{equation}
By~\eqref{eq:CC-differ} and~\eqref{eq:D(C)}, the equation
\[
\d_{\fc}(s)=0
\mod \CV^2_{r(s)-1}+\tilde\CV^2_{r(s)}
\]
for semilinear $s$ implies $d_h(\om(s,1)^\alpha)=0$.
\end{proof}

\begin{lemma}\label{lem:linear}
Any cocycle $s\in\CV^1$ with $r(s)\ge 1$ is semilinear.
\end{lemma}
\begin{proof}
If there is a nonlinear coefficient, we can fix the values of $x_i$, $v_I^{\alpha,A}$, $|I|<r(s)$,
such that some form $\om(s,0)^\alpha$ is not linear with respect to $u_I^{\alpha,A}$, $|I|=r(s)$,
for these values of the parameters. Proposition~\ref{prop:hor-triv} implies
\begin{equation}
\label{H1=0}
        H_h^1(J^\infty(\pi'))=0.
\end{equation}
By Lemma~\ref{lem:d(w)}~\eqref{lem:d(w)-1} and~\eqref{H1=0}, we have $\om(s,0)^\alpha=d_h(w')$
for some function dependent only on $u_I^{\alpha,A}$ with $|I|=r(s)-1$. Taking into
account~\eqref{d_h} and~\eqref{eq:tot-der}, we see that $\om(s,0)^\alpha$ is linear.
\end{proof}

\begin{lemma}\label{lem:F_1}
For each cocycle $s\in\CV^1$ there exists a cohomological cocycle with coefficients in $\CF_1$.
\end{lemma}
\begin{proof}
Assume that $r(s)\ge 2$. By Lemma~\ref{lem:linear}, Lemma~\ref{lem:d(w)}~\eqref{lem:d(w)-2},
and~\eqref{H1=0}, we have
\begin{equation}\label{eq:exact}
\om(s,1)^\alpha=d_h(w_\alpha)
\end{equation}
for some functions $w_\alpha\in\Ci(\Ji(\pi'))$  dependent also on the parameters $x_i$,
$v_I^{\alpha,A}$ with $|I|\le r(s)-2$. Moreover, we can assume that $w_\alpha$ depend only on
the coordinates $u_I^{\alpha,A}$ with $|I|\le r(s)-1$. Replace the variables $u_I^{\alpha,A}$
in $w_\alpha$ by $v_I^{\alpha,A}$ and denote the obtained functions by $\tilde w_\alpha$.
Consider the vector field
\[
s'=\sum_\alpha\tilde w_\alpha D_{v^\alpha}\in \CV^0_{r(s)-1}.
\]
Then from~\eqref{eq:exact} and~\eqref{eq:diff-mod1} we have $r(s-\d_{\fc}(s'))\le r(s)-1$.
Proceeding in this way, by induction on $r(s)$ one completes the proof.
\end{proof}

\begin{proof}[Proof of Theorem~\textup{\ref{thm:main}}]
Let us first prove~\eqref{0ac}. Let $q\in\CV$.
If $r(q)>0$ then $\d_\fc(q)\neq 0$ by~\eqref{eq:diff-mod1}.
If $r(q)=0$, i.e.,
\[
q=\sum_\alpha f\alpha(x_1,\dots,x_n,v^1,\dots,v^m)D_{v^\alpha},
\]
then $\d_\fc(q)\neq 0$ is proved by direct computation.

Statement~\eqref{0ac} means that it is sufficient to prove~\eqref{1ac} locally.

Let $s\in\CV^1$ be a cocycle. By Lemmas~\ref{lem:linear} and~\ref{lem:F_1}, there is
a cohomological cocycle
\[
s'=\sum_{i,\beta}f^\beta_i\,dx_i\otimes D_{v^\beta}\in\CV^1
\]
such that the functions $f^\beta_i$ belong to $\CF_1$ and are linear with respect to the
coordinates $v_I^{\alpha,A}$, $|I|=1$. {}{}From $d_h\om(s',0)=0$ we obtain
\[
f^\beta_i = f_{i,\beta}(x_1,\dots,x_n,v^1,\dots,v^m)+
\sum_{\gamma,A}f^{\gamma,A}_\beta(x_1,\dots,x_n,v^1,\dots,v^m)v_i^{\gamma,A}.
\]
By~\eqref{eq:CC-differ}, the equation $\d_{\fc}(s')=0$ is equivalent to
\begin{equation} \label{eq:mod2}
D_i(f^\beta_j)-\sum_{\gamma=1}^mv^{\beta,\gamma}_i\, f^\gamma_j
=D_j(f^\beta_i)-\sum_{\gamma=1}^mv^{\beta,\gamma}_j\,f^\gamma_i
\end{equation}
for all $\beta$, $i$ and $j$. Comparing equations~\eqref{eq:fce} and~\eqref{eq:mod2}, one shows
that the functions $f^{\gamma,A}_\beta$ vanish for $|A|>1$. Then a straightforward computation
implies that $s'=\d_{\fc}(F)$ for some $F\in\CV^0_0$.
\end{proof}

\section{Flat representations}\label{sec:flat-rep}

\subsection{Definition and the main property}\label{subs:fc-def}
Consider a formally integrable infinitely prolonged differential equation $\xi\colon\CE\to M$
with the Cartan distribution $\CC\Dr\subset\Dr(\CE)$, where $M$ is the space of independent
variables. A \emph{flat representation} of $\CE$ is given by a fibering $\pi\colon M\to N$ and
a flat connection $\CF$ in the composition bundle $\eta=\pi\circ\xi\colon\CE\to N$ such that
for each $X\in\Dr(N)$
\begin{equation}\label{ccd}
  \CF(X)\in \CC\Dr.
\end{equation}
Let $x_1,\dots,x_{n'}$ be local coordinates in $N$ and $y^1,\dots,y^m$ be coordinates along
the fibers of the bundle $\pi$. Then any flat representation is determined by the correspondence
\begin{equation}\label{eq:fc-loc}
\d x_i\mapsto D_{x_i}+\sum_{\alpha=1}^ma_i^\alpha D_{y^\alpha},\quad i=1,\dots,n',
\end{equation}
where $D_{x_i}$, $D_{y^\alpha}$ are the total derivatives on $\CE$ and the functions
$a_i^j\in\Ci(\CE)$ satisfy the condition
\begin{equation}\label{eq:fc-cond}
[D_{x_i}+\sum_{\alpha=1}^ma_i^\alpha D_{y^\alpha},D_{x_j}+\sum_{\alpha=1}^ma_j^\alpha D_{y^\alpha}]
=0,\quad 1\le i<j\le n'.
\end{equation}

Let $\CE_\fc$ be the equation of flat connections in the bundle $\pi$ and consider a morphism
$\vf\colon\CE\to\CE_\fc$ such that the diagram
\begin{diagram}[LaTeXeqno]\label{diag:fc}
\CE&              &\rto^{\vf} &              &\CE_\fc\\
     &\rdto_{\xi}& &\ldto&\\
     &              &M&              &\\
\end{diagram}
is commutative. Take a point $\theta\in\CE$. Then $\vf(\theta)$ is the infinite jet of some flat
connection $\nabla_\theta$ in the bundle $\pi$. This means that any tangent vector to $N$ at
the point $\eta(\theta)$ can be lifted by $\nabla$ to $M$ and then by the Cartan connection in
$\CE$ to a tangent vector to $\CE$ at the point $\theta$. Evidently, this procedure is independent
of the choice of $\nabla$ and we obtain a flat representation for~$\CE$.

The converse construction is also valid, i.e., any flat representation determines a morphism $\vf$
satisfying the commutativity condition~\eqref{diag:fc}. To establish this fact let us consider
a point $\theta\in \CE$ of the equation $\CE$ and a tangent vector $w$ to $N$ at the point
$\eta(\theta)$. Then we can construct a vector $\theta_{\CF}(w)\in T_{\xi(\theta)}M$ by
\begin{equation}\label{eq:pointwise}
\theta_{\CF}(w)= \xi_*\bigl|_\theta(\CF(w))\in T_{\xi(\theta)}M,
\end{equation}
and one has $\pi_*(\theta_{\CF}(w))=w$. Assume now that $\CE\sbs\Ji(\zeta)$ for some bundle
$\zeta$ over $M$ and the connection $\CF$ can be extended to some connection $\tilde\CF$
in the bundle $\pi\circ\zeta_\infty$ in such a way that $\tilde\CF$ also possesses
property~\eqref{ccd} ($\tilde\CF$ may not be flat). This can be always done. For a point
$\theta\in\CE$,
consider the points $x=\xi(\theta)\in M$, $x'=\pi(x)\in N$ and a neighborhood $\CU'\sbs N$ of
$x'$ in $N$. Set $\CU=\pi^{-1}(\CU')\sbs M$. Let $\theta=[f]_x^\infty$, $f\in\Gal(\zeta)$ and
$\CU_f=j_\infty(f)(\CU)\sbs\Ji(\zeta)$. Then, using~\eqref{eq:pointwise}, for any point
$\tilde{x}\in\CU$ and a vector field $X$ on $\CU'$ we set
\[
(\nabla_\theta X)_{\tilde{x}}=\tilde{\theta}_{\tilde{\CF}}(X_{\tilde{x}'}),
\]
where $\tilde{\theta}=j_\infty(f)(\tilde{x})$ and $\tilde{x}'=\pi(\tilde{x})$. It is
straightforward to check that $\nabla_\theta$ is a local connection in $\pi$, its infinite jet
at $x$ depends on $\theta$ only and lies in $\CE_\fc$. The map $\vf\colon\theta\mapsto
[\nabla_\theta]_x^\infty$ is a morphism of $\CE$ to $\CE_\fc$ and the corresponding flat
representation coincides with $\CF$. Thus we obtain
\begin{theorem}\label{thm:fc}
Let $\xi\colon\CE\to M$ be an infinitely prolonged equation and $\pi\colon M\to N$ be a fiber
bundle. Then flat representations of $\CE$ in the bundle $\pi$ are in one-to-one correspondence
with the morphisms of $\CE$ to the equation of flat connections for $\pi$ satisfying commutativity
condition~\eqref{diag:fc}.
\end{theorem}

\begin{remark}\label{rem:fc}
In local coordinates the construction of the morphism $\vf$ looks quite simple. Namely, let
a flat representation be given by the functions $a_i^\alpha$, see~\eqref{eq:fc-loc}. Then we set
$\vf^*(v_i^\alpha)=a_i^\alpha$ and take all differential consequences of these relations, i.e.,
\[
\vf^*(v_{i,\si}^\alpha)=D_\si(a_i^\alpha).
\]
Due to relations~\eqref{eq:fc-cond}, these formulas are consistent. Here $v_i^\alpha$ are the
coordinates introduced at the end of Subsection~\ref{subs:more-conn}.
\end{remark}

\subsection{Deformations of a flat representation}\label{subs:deform}
Integrable equations possess parametric families of flat representations.
Let us study this structure.

A smooth family of flat connections $\CF(\ep)$ in the bundle $\eta$ is called a \emph{deformation}
of a flat representation $\CF^0$ if $\CF(\ep)$ meets~\eqref{ccd} and $\CF(0)=\CF^0$. Below we also
consider \emph{formal deformations}, which are not smooth families, but formal power series
in~$\ep$.

Set $U=\bar U_{\CF^0}\in\Dr(\La^1(\CE))$. Consider the submodule $\CV_{\CE}\subset\CC\Dr$ of
$\eta$\h vertical Cartan fields and set $\CV_{\CE}^q=\mu^*(\Lambda^q(N))\otimes\CV\subset
\Dr(\Lambda_h^q(\CE))$.

Let $\CF(\ep)$ be a (formal) deformation of $\CF^0=\CF(0)$. {}From~\eqref{ccd} it follows that
$\bU_{\CF(\ep)}=U+U_1\ep+O(\ep^2)$, where $U_1\in\CV_\CE^1$. By Proposition~\ref{prop:curvature},
we have
\[
\fnij{U+U_1\ep+O(\ep^2)}{U+U_1\ep+O(\ep^2)}=2\fnij{U}{U_1}\ep+O(\ep^2)=0.
\]
Therefore, the \emph{infinitesimal part} $U_1\in \CV_{\CE}^1$ of any deformation satisfies the
identity
\begin{equation}\label{U1}
  \fnij{U}{U_1}=0.
\end{equation}

For any vector field $V\in \CV_{\CE}$ we can consider the formal deformation
\begin{equation}\label{ftr}
  \mathrm{e}^{\ep\ad V}(U)=
U+\fnij{V}{U}\ep+\frac12\fnij{V}{\fnij{V}{U}}\ep^2+\dots
\end{equation}
with the infinitesimal part
\begin{equation}\label{itr}
  \fnij{V}{U}.
\end{equation}
Since such deformations exist independently of the equation structure, it is natural to call
them, as well as their infinitesimal parts $\fnij{V}{U}$, $V\in\CV_\CE$, \emph{trivial}.

Formulas~\eqref{U1} and~\eqref{itr} prompt the following construction. It is easily seen that
$\d_{U}(\CV_{\CE}^q)\subset\CV_\CE^{q+1}$, and, therefore, we obtain a complex
\begin{equation}\label{eq:compl_gen}
0\to\CV_\CE^0=\CV_{\CE}
\xrightarrow{\d_{U}}\CV_{\CE}^1\xrightarrow{\d_{U}}\dots
\to\CV_\CE^q\xrightarrow{\d_{U}}\CV_\CE^{q+1}\to\dotsb
\end{equation}
By the construction and formulas~\eqref{U1}, \eqref{itr}, its $1$\h cocycles are infinitesimal
deformations, and exact $1$\h cocycles are trivial infinitesimal deformations. {}From the general
theory of deformations~\cite{gen-deform}, we obtain that obstructions to prolongation of infinitesimal deformations to
formal ones belong to higher cohomology groups of~\eqref{eq:compl_gen}.
\begin{remark}
Deformations of \emph{zero\h curvature representations} (finite\h dimensional linear coverings) of
PDE were studied analogously in~\cite{marvan02}. In fact, the \emph{gauge complex} of a
zero\h curvature representation introduced in~\cite{marvan92,marvan02} is a particular case of the
above construction, since according to Subsection~\ref{c-fr} any finite\h dimensional covering
determines a flat representations.
\end{remark}

\subsection{Symmetries and cocycles}\label{sc}
Let $S\in\Dr(\CE)$ be a symmetry of $\CE$. Since $S$ commutes with the fields of the form
$\CC(X)$, $X\in\Dr(M)$, the power series $\mathrm{e}^{\ep\ad S}(U)$ is a formal deformation of
$\CF$, and its infinitesimal part $\fnij{S}{U}\in\CV_\CE^1$ is a $1$\h cocycle of
complex~\eqref{eq:compl_gen}. Let us set $c_S=\fnij{U}{S}$. Note that  this cocycle is not
generally exact, since $S\notin \CV_\CE$.

\begin{theorem}
Symmetries $S$ with exact cocycles $c_S$ form a Lie subalgebra in $\sym\CE$.
\end{theorem}
\begin{proof}
Let $S_1,\,S_2$ be two symmetries such that
\begin{equation}\label{sivi}
\fnij{U}{S_1}=\fnij{U}{V_1},\quad \fnij{U}{S_2}=\fnij{U}{V_2},\quad V_1,\,V_2\in\CV_\CE.
\end{equation}
Using this and properties~\eqref{eq:skew}, \eqref{eq:Jacobi} of the Nijenhuis bracket, we easily
obtain
\begin{equation}\label{ssvv}
        c_{[S_1,S_2]}=\fnij{U}{[S_1,S_2]}=\\
        \fnij{U}{[V_1,V_2]+[V_1,S_2]+[S_1,V_2]}.
\end{equation}
Note that $[\CV_\CE,\CV_\CE]\subset\CV_\CE$ and for each symmetry $S$ one has $[S,\CV_\CE]\subset
\CV_\CE$. Therefore, $[V_1,V_2]+[V_1,S_2]+[S_1,V_2]\in\CV_\CE$ and the cocycle $c_{[S_1,S_2]}$ is
exact.
\end{proof}

\begin{theorem}\label{sym_gen}
Suppose that the morphism $\varphi\colon\CE\to\CE_{\fc}$ corresponding to the flat representation
$\CF$ is an embedding such that $\varphi(\CE)$ is a submanifold of $\CE_{\fc}$. Then the cocycle
$c_S$ is exact if and only if there is a symmetry $S'\in\Dr(\CE_{\fc})$ of $\CE_{\fc}$ such that
\begin{equation}
\label{restrict}
\varphi_*(S)=S'\bigl|_{\varphi(\CE)}.
\end{equation}
\end{theorem}
\begin{remark}
Recall that $\CE$ is a submanifold of some infinite jet space $\Ji(\zeta)$. It is well
known~\cite{KV-book} that under rather weak regularity assumptions each symmetry of $\CE$ is the
restriction of some symmetry of $\Ji(\zeta)$. The above theorem provides a criterion for the
similar property in the situation when $\CE_{\fc}$ is considered instead of $\Ji(\zeta)$.
In the next section we give examples of both types of symmetries: with exact and non-exact $c_S$.
Therefore, generally not every symmetry is the restriction of some symmetry of $\CE_{\fc}$.
\end{remark}
\begin{proof}
First suppose that~\eqref{restrict} holds.
By Theorem~\ref{thm:main}, there is $V\in\CV$ such that
\begin{equation}\label{s'v}
        \fnij{U}{S'}=\fnij{U}{V}.
\end{equation}
Since $V$ belongs to the Cartan distribution and $\varphi$ is a morphism of differential equations,
the vector field $V$ is tangent to the submanifold $\varphi(\CE)$. Denote by $V'\in\CV_{\CE}$
the preimage $\varphi_*^{-1}(V)$. Then from~\eqref{restrict} and~\eqref{s'v} it follows that
\begin{equation}\label{sv'}
        c_S=\fnij{U}{V'},
\end{equation}
i.e., the cocycle $c_S$ is exact.

Conversely, if $c_S$ is exact then~\eqref{sv'} holds for some $V'\in\CV_{\CE}$. Extend
$\varphi_*(V')\in \Dr(\varphi(\CE))$ to a vector field $V\in\CV$. Then the symmetry $S'$ of
$\CE_{\fc}$ with the generating section $\fnij{U}{V}$ satisfies~\eqref{restrict}.
\end{proof}

\subsection{A bracket on the equation of flat connections}\label{subs:Lie-str}
The identity mapping
\[
\id\colon\CE_{\fc}\to\CE_{\fc}
\]
is a flat representation, and complex~\eqref{eq:CC-Vcompl} coincides with~\eqref{eq:compl_gen} for
this flat representation. {}From the results of Section~\ref{subs:result} it follows that
$\sym \CE_{\fc}$ is isomorphic to $\CV$. Therefore, the commutator of symmetries induces a bracket
on $\CV$. By formulas~\eqref{sivi} and~\eqref{ssvv}, it is given by
\begin{equation}\label{bracket}
\{V_1,V_2\}=[V_1,V_2]+[V_1,S_2]+[S_1,V_2],\quad V_1,\,V_2\in\CV,
\end{equation}
where $S_i$ is the symmetry with the generating section $\fnij{\bU_{\fc}}{V_i}$. Let us compute
this bracket in coordinates.

In adapted coordinates, any element of $\CV$ is of the form
\[V_f=\sum_{i=1}^nf^iD_{v^i},\quad f^i\in\Ci(\CE_\fc).
\]
Hence, locally, the bracket $\{V_f,V_g\}$ generates a bracket $\{f,g\}$ on vector\h functions
defined by $\{V_f,V_g\}=V_{\{f,g\}}$. Using~\eqref{bracket}, we compute $\{f,g\}$.

Let
\[
S_f=\sum_{|I|>0} S_I^{\alpha,A}\d v_I^{\alpha,A}
\]
be the symmetry of $\CE_\fc$ with the generating section $\fnij{\bU_{\fc}}{V_f}$. Then by direct
computations we obtain that
\begin{equation}\label{eq:v-comp}
S_I^{\alpha,A}=D_{v^A}(S_I^{\alpha,\varnothing}),
\end{equation}
where $D_{v^A}$ is the composition of the total derivatives $D_{v^\beta}$ corresponding to the
multi\h index $A$, while
\begin{equation}\label{eq:x-comp}
S_{Ii}^{\alpha,\varnothing}=D_i(S_I^{\alpha,\varnothing})
+\sum_{\beta=1}^mv_I^{\alpha,\beta}\vf_i^\beta,
\end{equation}
where $\vf_i^\beta= S_i^{\beta,\varnothing}$, and, as we already know (see~\eqref{eq:spec-coord}),
\[
\vf_i^\beta=D_i(f^\beta)-\sum_{\gamma=1}^mv_i^{\beta,\gamma}f^\gamma.
\]

Substituting these expressions to~\eqref{bracket}, we obtain the following
\begin{proposition}\label{prop:Lie-str}
The commutator of symmetries determines the following bracket on the module of smooth
$m$\h dimensional vector\h functions on~$\CE$
\begin{equation}
\{f,g\}^\alpha=S_fg^\alpha-S_gf^\alpha+V_fg^\alpha-V_gf^\alpha.
\end{equation}
For any vector\h function $h=(h^1,\dots,h^m)$ the coefficients of the vector field
$S_h=\sum S_I^{\alpha,A}\d v_I^{\alpha,A}$ in special coordinates $v_I^{\alpha,A}$ can be
computed by the formulas
\[
S_I^{\alpha,A}=D_{v^A}(S_I^{\alpha,\varnothing}),\qquad
S_{Ii}^{\alpha,\varnothing}=D_i(S_I^{\alpha,\varnothing})
+\sum_{\beta=1}^mv_I^{\alpha,\beta}\zeta_i^\beta,
\]
where $\zeta_i^\beta= S_i^{\beta,\varnothing}$ and $\zeta_i^\beta=D_i(h^\beta)-
\sum_{\gamma=1}^mv_i^{\beta,\gamma}h^\gamma$\textup{,} while $V_h=
\sum_{\alpha}h^\alpha D_{v^\alpha}$.
\end{proposition}

\section{Examples of flat representations}\label{efr}

\subsection{Coverings as flat representations and lifting of symmetries}\label{c-fr}

Consider an infinitely prolonged differential equation $\xi\colon\CE\to M$
with the Cartan connection
\[
\CC\colon\Dr(M)\to\Dr(\CE).
\]
Let $\pi\colon E\to M$ be a fiber bundle of finite rank and consider its
pullback with respect to~$\xi$
\begin{equation}\label{cov}
\begin{CD}
        \CE'@>\xi'>> E\\
        @ V\xi^*(\pi)VV @ VV\pi V\\
        \CE @ >\xi>> M
\end{CD}
\end{equation}
Suppose that the bundle $\tau=\xi^*(\pi)\colon\CE'\to\CE$ is endowed with a covering structure
given by a flat connection $\CF\colon\Dr(M)\to\Dr(\CE')$, $\tau_*\circ\CF=\CC$. This structure
determines a flat representation not of the equation~$\CE$ itself, but of some trivial extension
of~$\CE$.

Namely, consider the distribution $\CC D_1=\tau_*^{-1}(\CC\Dr_\CE)\subset\Dr(\CE')$. Clearly, the
differential $\xi'_*$ projects $\CC\Dr_1$ isomorphically onto $\Dr(E)$. Denote by $\CC_1\colon
\Dr(E)\to\Dr(\CE')$ the arising connection in the bundle $\xi'\colon\CE'\to E$. It is easily seen
that this connection is flat, and, moreover, the pair $(\CE',\CC_1)$ is isomorphic to an
infinitely prolonged equation with $E$ as the space of independent variables and $\CC_1$ as the
Cartan connection.

\begin{remark}
In coordinates, this construction looks as follows. Locally the bundle $\tau$ is trivial
$\tau\colon M\times W\to M$. Let $x_i$, $w_k$ be coordinates in $M$ and $W$ respectively.
Then the equation $(\CE',\CC_1)$ is obtained from $(\CE,\CC)$ as follows. We assume that the
dependent variables $u^j$ in $\CE$ are functions of not only $x_i$, but also $w_k$, and add
the equations $\d u^j/\d w_k=~0$.
\end{remark}

The bundle $\pi\colon E\to M$ and the connection $\CF$ constitute a flat representation of the
equation $(\CE',\CC_1)$.

Consider complex~\eqref{eq:compl_gen} corresponding to this flat representation. By construction,
the module $\CV_{\CE}$ consists of $\tau$\h vertical vector fields on~$\CE'$
and~\eqref{eq:compl_gen} in this case is the horizontal de~Rham complex of the covering equation
$(\CE',\CF)$ with coefficients in $\CV_{\CE}$ (see Section~\ref{rem:hor-coeff}). For arbitrary
coverings this complex was studied in~\cite{backlund} in relation with parametric families of
B\"acklund transformations.

By construction~\eqref{cov}, each $\xi$\h vertical vector field $Y\in\Dr(\CE)$ is canonically lifted
to $\CE'$ as a $\xi'$\h vertical vector field denoted by $\hat{Y}\in\Dr(\CE')$. Evidently, for each
symmetry $S$ of $\CE$ the vector field $\hat{S}$ is a symmetry of the equation $(\CE',\CC_1)$.

\begin{theorem}\label{lifts}
Let $S$ be a symmetry of $\CE$ and consider complex~\eqref{eq:compl_gen} corresponding to the flat
representation $\CF$ of $(\CE',\CC_1)$. The cocycle $c_{\hat S}$ is exact if and only if there is
a symmetry $S'\in\Dr(\CE')$ of the covering equation $(\CE',\CF)$ such that $\tau_*(S')=S$.
\end{theorem}
\begin{proof}
If $c_{\hat S}=\fnij{\bar U_{\CF}}{\hat S}=\fnij{\bar U_{\CF}}{V}$ for some $V\in \CV_{\CE}$ then
\[
\fnij{\bar U_{\CF}}{\hat S-V}=0,
\]
i.e., the vector field $S'=\hat S-V$ is a symmetry of $(\CE',\CF)$. Clearly, $\tau_*(S')=S$.

Conversely, if $S'$ is a symmetry of $(\CE',\CF)$ and $\tau_*(S')=S$ then
$c_{\hat S}=\fnij{\bar U_{\CF}}{V}$ for $V=\hat S-S'\in\CV_\CE$.
\end{proof}

\subsection{Flat representations and symmetries of the KdV equation}\label{sec:flat-KdV}
Let us illustrate the above construction by the following classical example of covering.

Consider the well-known system
\begin{align}
\label{kdva}
        v_x&=\la+u+v^2,\\
\label{kdvb}
        v_t&=u_{xx}+2u^2-2\la u-4\la^2+2u_xv+v^2(2u-4\la),\,\ \la\in\mathbb{R},
\end{align}
associated with the Miura transformation.
Its compatibility condition is equivalent to the KdV equation
\begin{equation}\label{kdv}
u_t=u_{xxx}+6uu_x.
\end{equation}
The infinite prolongation $\CE$ of~\eqref{kdv} has the natural coordinates
\begin{equation}\label{coorkdv}
        x,\ t,\ u_k=\frac{\d^k u}{\d x^k},\quad k=0,1,2,\dots
\end{equation}
The total derivatives restricted to $\CE$ are written in these coordinates as follows
\begin{align}\label{Dx}
        D_x&={\d x}+\sum_{k\ge 0}u_{k+1}\,{\d u_k},\\
\label{Dt}
        D_t&={\d t}+\sum_{k\ge 0}D_x^k(u_3+6uu_1)\,{\d u_k}.
\end{align}
Evidently, functions~\eqref{coorkdv} and $v$ form a system of coordinates for the infinite
prolongation $\CE'=\CE\times\BBR$ of system~\eqref{kdva},~\eqref{kdvb}. In these coordinates, the
total derivatives restricted to $\CE'$ are $\tilde D_x=D_x+A{\d}/{\d v}$, $\tilde D_t=
D_x+B{\d}/{\d v}$, where $A,\,B\in\Ci(\CE')$ are the right-hand sides of~\eqref{kdva}
and~\eqref{kdvb} respectively. We have the trivial bundle  $\tau\colon\CE'\to\CE$, $(x,t,u_k,v)
\mapsto (x,t,u_k)$, equipped with the covering structure $\tau_*(\tilde D_x)=D_x$,
$\tau_*(\tilde D_t)=D_t$.

By the construction of the previous subsection, this covering determines the following flat
representation of the KdV equation.

We assume that $u$ in~\eqref{kdv} is a function of $x,\,t,\,v$ and add the condition $\d u/\d v=0$.
The infinite prolongation of this trivially extended KdV equation is the manifold $\CE\times\BBR$
with coordinates~\eqref{coorkdv} and $v$. The Cartan distribution $\CC\Dr_1$ on $\CE\times\BBR$ is
$3$\h dimensional. The total derivatives with respect to $x,\,t$ are given by~\eqref{Dx}
and~\eqref{Dt}, while the total derivative with respect to $v$ is just $\d v$. Then a flat
representations dependent on parameter $\la$ is given by the flat connection
\[
{\d x}\mapsto D_x+A{\d v}\,\in\CC D_1,\quad {\d t}\mapsto D_t+B{\d v}\,\in\CC\Dr
\]
in the bundle $\CE\times\BBR\to\BBR^2$.

The infinitesimal part of this parametric family is
\begin{equation*}
        \frac{\d A}{\d \la}\,dx+\frac{\d B}{\d \la}\,dt=
        dx\otimes\d v-(2u+8\la+4v^2)\,dt\otimes\d v.
\end{equation*}
This $1$\h cocycle is not exact, which reflects the well-known fact that the parameter here is
essential.

According to Theorem~\ref{lifts}, the cocycle $c_S$ corresponding to a symmetry $S$ of~\eqref{kdv}
is exact if and only if the symmetry is lifted to a symmetry of system~\eqref{kdva},~\eqref{kdvb}.
It is well known which symmetries of~\eqref{kdv} can be lifted. Namely, all $(x,t)$\h independent
symmetries (including the symmetries corresponding to the higher KdV equations) are lifted, while
the Galilean symmetry is not. The scaling symmetry can be lifted only in the case $\la=0$.

\begin{remark}
For an equation $\CE$ in two independent variables existence of a nontrivial covering is an
important indication of integrability and often leads to B\"acklund transformations and the inverse
scattering method~\cite{KV-trends}. For non\h overdetermined equations in more than two
independent variables there are no nontrivial finite\h dimensional coverings (this was proved
in~\cite{marvan92} for linear coverings and in \cite{igonin} for arbitrary coverings). The notion
of flat representation may serve as a good substitute in this case, see a typical example in the
next subsection.
\end{remark}

\subsection{Flat representations of the self-dual Yang--Mills equations}
Consider four smooth vector\h functions $A_i\colon\mathbb{R}^4\to\g$, $i=1,2,3,4$, where $\g$ is a
finite\h dimensional Lie algebra. The self-dual Yang--Mills equations are
\begin{equation} \label{sdym}
  [{\d x_1}+A_1+\la({\d x_3}+A_3),\,{\d x_2}+A_2+\la({\d x_4}+A_4)]=0,
\end{equation}
which must hold for all values of the parameter $\la$ (see, for example,~\cite{ward}). Here $x_i$
are coordinates in $\mathbb{R}^4$. Denote by $\CE$ the infinite prolongation of~\eqref{sdym} and
by $\CC$ the Cartan connection in the bundle $\xi\colon\CE\to\mathbb{R}^4$.

Fix an action $\sigma\colon \g\to\Dr(W)$ of $\g$ on a finite\h dimensional manifold $W$. Similarly
to Subsection~\ref{c-fr}, we construct a trivial extension $\CE'=\CE\times W$ of $\CE$ as follows.
Consider the flat connection
\begin{equation*}
  \CC_1=\CC\oplus\idt\colon\Dr(\BBR^4\times W)\to\Dr(\CE')
\end{equation*}
in the bundle $\xi_1=\xi\times\idt\colon\CE'\to \BBR^4\times W$.

Evidently, the pair $(\CE',\CC_1)$ is isomorphic to an infinitely prolonged equation with
$\BBR^4\times W$ as the space of independent variables and $\CC_1$ as the Cartan connection.
Consider the bundle
\begin{gather*}
\pi\colon\BBR^4\times W\to \BBR^2,\,(a,w)\mapsto (x_1,x_2)\in\BBR^2,\\
a=(x_1,x_2,x_3,x_4)\in\BBR^4,\ w\in W.
\end{gather*}
The bundle $\pi$ and the flat connection $\CF$ in the bundle $\pi\circ\xi_1$ given by
\begin{equation}\label{frsdym}
\CF({\d x_i})=\CC_1\bigl({\d x_i}+\sigma(A_i)+\la({\d x_{i+2}}+\sigma(A_{i+2}))\bigl),\quad i=1,2,
\end{equation}
constitute a $\la$\h dependent family of flat representations for the (trivially extended)
self-dual Yang--Mills equations~\eqref{sdym}.

The $1$\h cocycle $\CC_1\bigl({\d x_3}+\sigma(A_3)\bigl)\,dx^1+\CC_1\bigl({\d x_4}+
\sigma(A_4)\bigl)\,dx^2$ corresponding to the family of flat representations~\eqref{frsdym} is not
exact, which says that the parameter in~\eqref{frsdym} is essential.

Let $H\colon\CE\to\g$ be a smooth function; then the formula
\begin{equation}\label{ggsym}
        G_H(A_i)=D_{x_i}(H)-[H,A_i]
\end{equation}
defines a higher symmetry of $\CE$. Since in the case when $H$ depends on $x_i$ only (does not
depend on $A_i$ and their derivatives) this is a classical gauge symmetry, it is natural to
call~\eqref{ggsym} the \emph{generalized gauge symmetry} corresponding to~$H$.

By the definition in Section~\ref{sc}, the $1$\h cocycle corresponding to the symmetry $G_H$ is
\begin{multline}\label{ugh}
        \fnij{\bar U_{\CF}}{G_H}=\sum_{i=1,\,2}[D_{x_i}+\sigma(A_i)+
\la(D_{x_{i+2}}+\sigma(A_{i+2})),G_H]\,dx^i=\\
  -\sigma\bigl(G_H(A_1+\la A_3)\bigl)\,dx^1-\sigma\bigl(G_H(A_2+\la A_4)\bigl)\,dx^2=
  -\fnij{\bar U_{\CF}}{\sigma(H)}
\end{multline}
The last equality in~\eqref{ugh} says that this cocycle is exact.

On the other hand, it can be shown that cocycles corresponding to classical conformal symmetries
of~\eqref{sdym} (see, for example,~\cite{popov}) are not exact.

\begin{remark}
Some other symmetries of~\eqref{sdym} (see~\cite{popov} and references therein) have been
extensively studied, but they are \emph{nonlocal} and, therefore, lie out of scope of the present
article.
\end{remark}

\begin{remark}
Various integrable reductions of the self-dual Yang--Mills equations possess similar
to~\eqref{sdym} Lax pairs (see~\cite{soliton,ward} and references therein) and, therefore, admit
flat representations.
\end{remark}

\section*{Acknowledgements}
The third author is grateful for hospitality to the University of Twente, where this paper was
started and finished.


\begin{thebibliography}{99}

\bibitem{soliton}
M.J.~Ablowitz and P.A.~Clarkson, \emph{Solitons, nonlinear evolution equations and inverse
scattering}, Cambridge Univ. Press, Cambridge, 1991.

\bibitem{anderson}
I.M.~Anderson and C.G.~Torre,
Classification of local generalized symmetries for the vacuum
Einstein equations, \emph{Comm. Math. Phys.} \textbf{176} (1996) no.~3, 479--539.

\bibitem{KV-book}
A.V.~Bocharov, V.N.~Chetverikov, S.V.~Duzhin, N.G.~Khor{\cprime}kova, I.S.~Krasil{\cprime}shchik,
A.V.~Samokhin, Yu.N.~Torkhov, A.M.~Verbovetsky, and  A.M.~Vinogradov,
\emph{Symmetries and Conservation Laws for Differential Equations of Mathematical Physics},
Amer. Math. Soc., Providence, RI, 1999. Edited and with a preface by Krasil{\cprime}shchik
and Vinogradov.

\bibitem{gen-deform}
M.~Gerstenhaber and S.~Schack,
Algebraic cohomology and deformation theory, in: M.~Hazewinkel, (Ed.),
\emph{Deformation theory of algebras and structures and applications},
Kluwer Acad. Publ., Dordrecht, 1988, pp.~11--264.

\bibitem{igonin}
S.~Igonin, Horizontal cohomology with coefficients and nonlinear
zero-curvature representations, \emph{Russian Math. Surveys}, to appear.

\bibitem{backlund}
S.~Igonin and J.~Krasil{\cprime}shchik,
On one-parametric families of B\"{a}cklund transformations,
in: Tohru Morimoto, Hajime Sato, and Keizo Yamaguchi, (Eds.)
\emph{Lie Groups\textup{,} Geometric Structures and Differential Equations\textup{:}
One Hundred Years after Sophus Lie},
Advanced Studies in Pure Mathematics, \textbf{37} (2002), p. 99--114; 
\url{arXiv.org/nlin.SI/0010040}

\bibitem{K-newinv}
I.S.~Krasil{\cprime}shchik,
Some new cohomological invariants for nonlinear differential equations,
\emph{Differential Geom. Appl.} \textbf{2} (1992) no.~4, 307--350.

\bibitem{Kras-flat}
I.S.~Krasil{\cprime}shchik,
Algebras with flat connections and symmetries of differential equations,
in Lie groups and Lie algebras, in: B.~Komrakov, I. Krasil{cprime}shchik, and A.~Sossinsky,
(Eds.), \emph{Lie groups and Lie algebras}, Kluwer Acad. Publ., Dordrecht, 1998, pp.~407--424.

\bibitem{KK-book}
I.S.~Krasil{\cprime}shchik and P.H.M.~Kersten,
\emph{Symmetries and recursion operators for classical and supersymmetric differential equations},
Kluwer Acad. Publ., Dordrecht, 2000.

\bibitem{KrasVer}
I.S.~Krasil{\cprime}shchik and A.M.~Verbovetsky,
\emph{Homological methods in equations of mathematical physics},
Open Education and Sciences, Opava (Czech Republic), 1998;
\url{arXiv.org/math.DG/9808130}

\bibitem{KV-trends}
I.S.~Krasil{\cprime}shchik and A.M. Vinogradov,
Nonlocal trends in the geometry of differential equations: symmetries,
conservation laws, and B\"{a}cklund transformations,
\emph{Acta Appl. Math.} \textbf{15} (1989) no.~1-2, 161--209.

\bibitem{marvan92}
M.~Marvan,
On zero-curvature of partial differential equations,
Proc. Conf. on Diff. Geom. and Its Appl.,
Opava (Czech Republic), 1992, p.~103--122; 
\url{www.emis.de/proceedings/5ICDGA/}

\bibitem{marvan02}
M.~Marvan,
On the horizontal gauge cohomology and
non-removability of the spectral parameter,
\emph{Acta Appl. Math.} \textbf{72} (2002) 51--65;
\url{diffiety.org/preprint/2001/02_01abs.htm}

\bibitem{popov}
A.D.~Popov,
Self-dual Yang--Mills: symmetries and moduli space,
\emph{Rev. Math. Phys.}
\textbf{11} (1999) 1091--1149; 
\url{arxiv.org/hep-th/9803183}

\bibitem{ward}
R.~Ward, Integrable systems and twistors,
in: \emph{Integrable systems \textup{(}Oxford, 1997\textup{)}}, 121--134.
Oxford Univ. Press, New York, 1999.
\end{thebibliography}
\end{document}